%% file: thread-arxiv.tex
\documentclass[fleqn]{article}
\usepackage[totalwidth=5.75in,totalheight=9.25in]{geometry}

\usepackage[mathscr]{euscript}
\usepackage[pdftex]{hyperref}
\usepackage{amsfonts}
\usepackage{amsmath}
\usepackage{amssymb}
\usepackage{amsthm}
\usepackage{bm}
\usepackage{caption}
\usepackage{graphicx}
\usepackage{lineno}
\usepackage{mathrsfs}
\usepackage{mathtools}
\usepackage{pgfgantt}
\usepackage{pgfplots}
\usepackage{tikz}
\usepackage{titlesec}
\usepackage{todonotes}
\usepackage{xcolor}
\usepackage{enumitem}
\usepackage{todonotes}

\definecolor{refcolor}{RGB}{42,93,176}
\hypersetup{colorlinks=true,linkcolor=refcolor,citecolor=refcolor}

\titleformat*{\section}{\large\bfseries}
\titlespacing*\section{0pt}{10pt plus 2pt minus 1pt}{6pt plus 0pt minus 0pt}
\newtheorem{definition}{Definition}

\usetikzlibrary{calc,patterns,angles,quotes}
\pgfplotsset{compat=newest}
\pgfplotsset{plot coordinates/math parser=false}
\newlength{\fwidth}

\newcommand{\VarTay}{{VT}}
\newcommand{\GauLegFour}{{GL4}}
\newcommand{\GauLegSix}{{GL6}}
\newcommand{\RadThree}{{R3}}
\newcommand{\RadFive}{{R5}}
\newcommand{\trbdf}{{TRB}}
\newcommand{\CN}{{TR}}
\newcommand{\BE}{{BE}}

\newcommand{\FF}{\bm F}
\newcommand{\ff}{\bm f}

\newcommand{\qq}{\bm q}
\newcommand{\xx}{\bm x}
\newcommand{\yy}{\bm y}

\begin{document}\allowdisplaybreaks\mathtoolsset{showonlyrefs}

\thispagestyle{empty}
\setlist{nosep}

\begin{raggedright}
{\LARGE\bf
\input{0title}
}\\[15pt]
Kevin R.\ Green$^{1,3}$, George W.\ Patrick$^{2,3}$, and Raymond J. Spiter$^{1,3}$\\[5pt]
$^1$Department of Computer Science, University of Saskatchewan,
  Saskatoon SK S7N5C9 Canada\\
$^2$Department of Mathematics and Statistics, University of Saskatchewan,
  Saskatoon SK S7N5E6 Canada\\
$^3$Acknowledge support from the Natural Sciences and Engineering Research
  Council of Canada (NSERC)\\[5pt]
The authors thank the anonymous referees for their insightful comments.\\[5pt]
\today\\[10pt]
{\color{lightgray}\rule{\textwidth}{.5pt}}
\parbox{\textwidth}{
  \vspace*{3pt}{\noindent\bf Abstract\\[4pt]}
  \input{0abstract.tex}

}
\\[4pt]
{\color{lightgray}\rule{\textwidth}{.5pt}}\\[6pt]
\end{raggedright}
\vspace*{1pt}

\input{0section-1-introduction.tex}
\input{0section-2-continuation.tex}
\input{0section-3-double-pendulum.tex}
\input{0section-4-discussion.tex}

\footnotesize
\bibliography{references}

\end{document}

%% file: 0title.tex
On theoretical upper limits for valid timesteps of implicit ODE methods

%% file: 0abstract.tex
Implicit methods for the numerical solution of initial-value problems
may admit multiple solutions at any given time step. Accordingly,
their nonlinear solvers may converge to any of these solutions. Below
a critical timestep, exactly one of the solutions (the consistent
solution) occurs on a solution branch (the principal branch) that can
be continuously and monotonically continued back to zero timestep.

\smallskip

Standard step-size control can promote convergence to consistent solutions by
adjusting the timestep to maintain an error estimate below a given
tolerance. However, simulations for symplectic systems or large physical
systems are often run with constant timesteps and are thus more susceptible to
convergence to inconsistent solutions. Because simulations cannot be
reliably continued from inconsistent solutions, the critical
timestep is a theoretical upper bound for valid timesteps.


%% file: 0section-1-introduction.tex
\section{Introduction}
\label{sec:introduction}

Many mathematical models take the form of a system of ordinary differential
equations (ODEs) for a vector of unknowns $\qq(t) \in \mathbb{R}^m$ subject to
boundary data:
\begin{alignat*}{2}
  \dot\qq(t) &= {\ff}(\qq(t)),    &&\quad t \in (t_0,t_f),\\
           0 &= g_i(\qq(\tau_i)), &&\quad \tau_i \in \{t_0,t_f\}, \ i=1,2,\ldots, m.
\end{alignat*}
Standard transformations reduce systems that are higher order,
non-autonomous, or subject to interior-point data to this first-order
autonomous form with boundary data at the cost of increased system
dimension~\cite{Ascher1995}. When $\tau_i \equiv \tau_0$,
$i=1,2,\ldots, m$, we have an initial-value problem~(IVP); otherwise,
we have a two-point boundary-value problem~(BVP).

In practice, numerical methods for the solution for such ODEs involve
successive approximations at successive timesteps and are either
\emph{implicit} or \emph{explicit}. Implicit methods typically involve
the iterative solution, at each time step, of systems of nonlinear
algebraic equations, generally a theoretically infinite process with a
potentially non-unique or non-existent result.  Explicit methods, in
contrast, can be implemented directly, generally a theoretically
finite process with a unique result.  The iterative solution process
of an implicit method can incur a significant run-time cost, but the
use of such methods may result in greater overall efficiency or
fidelity. For example, the increase in the timestep afforded by an
implicit method when solving stiff ODEs typically offsets the
increased cost per step. Also, when integrating a Hamiltonian system,
an implicit method may be used to arrange that the simulation itself
preserves energy or is
symplectic~\cite{HairerE-LubichC-WannerG-2006-1,
  LeimkuhlerB-ReichS-2004-1}.

The existence and uniqueness theory for IVPs is much more decisive
than for BVPs. IVPs have unique solutions under mild assumptions that
are typically satisfied in practice, whereas BVPs may have from zero
to uncountably many solutions. However, when an implicit method is
involved in approximating the solution of an IVP, the possibility
emerges of divergence or convergence to one of multiple
solutions. Convergence to spurious solutions is well recognized
(\cite{Deuflhard2011,Stuart1998} and references therein),
particularly in the numerical solution of
BVPs~(\cite{Ascher1995,Murdoch1992,Schreiber1983,Stephens1981} and
references therein). Less attention is typically given to the context
of solving IVPs (\cite{Griffiths1992,Humphries1993,Iserles1991} and
references therein), where there is theoretically a unique solution
and where the presence of a ``good'' initial guess is taken for
granted.  In this article, we consider the context of an implicit IVP
method that has multiple solutions at a given timestep and how to
choose from among them, as opposed to a qualitatively incorrect
numerical solution of an IVP or BVP.


Standard methods for error estimation and control via timestep selection tend
to adjust timesteps such that, in practice, any ambiguity arising from multiple
solutions is avoided, but they are not usually specifically designed to do
so. But there are two specific scenarios in which constant timesteps are often
used in practice: simulations of symplectic systems using a symplectic
method~\cite{HairerE-LubichC-WannerG-2006-1} and simulations of large physical
systems. Efforts toward adaptive symplectic methods have been made, but they
tend to be specialized and require a significant amount of user
judgment~\cite{Hairer1997,Reich1999,Blanes2005,Hairer2005}. Software packages
for the simulation of large physical systems, especially on distributed
architectures~\cite{Phillips2005,Pitt-Francis2008,Cantwell2015,Juno2018}, often
use constant timesteps, because of the large relative expense of estimating the
error and changing the timestep. However, at constant timestep, a simulation
may unexpectedly encounter a time-localized region of complex dynamics, and
convergence can easily be construed as nominal when in fact it is not. A large
number of independent simulations, e.g., explorations of a parameter space, may
not be feasible or efficient at small constant timestep. It may not be easy to
automate the detection of anomalous behaviour in the absence of convergence
failure, and a manual inspection may not be feasible.



Numerical methods for solution of ODEs can also be classified
according to how many past steps are stored and used in computing the
next step. In general, the next step can be a function of $k$ past
steps, leading to \textit{multi-step methods}. If only the current
state is stored, then the method is \textit{one-step}. Although a
multistep method may be regarded as a one-step method on a Cartesian
product of state spaces~\cite{KirchgraberU-1986-1}, the theory of
multistep methods is complicated by the possibility of uncontrolled
growth of the error in the past states~\cite{IserlesA-2009-1}.  This
is not the focus of this article; we consider only one-step methods.

To formalize, the time-advanced state $\bm q^{n+1} \approx \qq(t_{n+1})$ after
one step of a one-step implicit method is obtained from the given current
state $\bm q^n \approx \qq(t_n)$ by solving a generally nonlinear equation of
the form
\begin{equation}\label{eq:general-time-integration1}
  \bm F(\bm y;\bm x,h)=\bm 0, \qquad
  \bm y=\bm q^{n+1},  \qquad
  \bm x= \bm q^n,
\end{equation}
where $h$ is the given timestep. For example, the backward
Euler method has $\FF(\yy;\xx,h) = \yy - \xx - h\,\ff(\yy)$.

We assume, for all $\bm x$, that
\begin{equation*}
  \bm F(\bm x;\bm x,0)=\bm 0,\qquad
  \bm F_{\bm y}(\bm x;\bm x,0)=\bm 1,\qquad
  \bm F_{h}(\bm x;\bm x,0)=-\bm f(\bm x),
\end{equation*}
where $\bm1$ is the identity matrix and subscripts of $\bm F$ denote partial
derivatives. Then, by the implicit function theorem, for any fixed $\bm x$,
there is a unique smooth solution $\bm y(h)$ defined for sufficiently small~$h$,
and $\bm y=\bm x+h\bm f(\bm x)+O(h^2)$, as is required for consistency.

Two solutions $\bm y_1(h)$ and $\bm y_2(h)$, defined on open intervals
containing $h=0$ and satisfying the condition that
$\bm F_{\bm y}\bigl(\bm y_i(h);\bm x,h\bigr)$ is nonsingular, are
equal on the intersection of their domains (the set on which
$\bm y_1(h)=\bm y_2(h)$ is nonempty, closed, and open by the implicit
function theorem, and the intersection of intervals is
connected). Therefore, there is a maximal such solution, which we call the
\emph{principal solution branch}, and there is generally a critical
timestep~$h_c$ after which the principal branch ceases to exist. The
condition that the solution set of $F(\bm y;\bm x,h)=\bm 0$ is smooth
in the space $\{\bm y,h\}$ is weaker than the condition that it
defines $\bm y$ as a function of $h$. Solutions may be continuously
connected after a fold bifurcation, for example, where the solution
manifold turns backwards from the direction of increasing
timestep~\cite{GuckenheimerJHolmesP-1983.1,
  KuznetsovYA-2004-1}. Continuing through such a bifurcation leads to
multiple solutions at smaller timesteps than the critical one at which
the bifurcation occurs. These solutions co-exist with the solutions
from the principal branch, and they may persist even as the timestep
approaches zero.

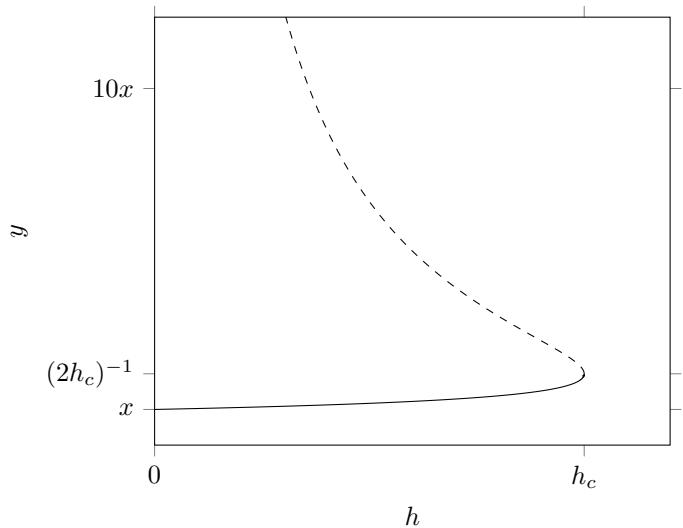
\begin{figure}[htbp]
  \centering
  \begin{tikzpicture}
    \begin{axis}[
            xmin=0.0,xmax=0.3,
            ymin=0.0,ymax=12,
            xtick align=outside,
            ytick align=outside,
            ylabel={$y$},
            xlabel={$h$},
            xtick={0,0.25},
            xticklabels={$0$,$h_c$},
            ytick={1,2,10},
            yticklabels={$x$,$(2h_c)^{-1}$,$10x$}
            ]
            \addplot [domain=1:2,samples=1000]({1/x-1/x^2},{x});
            \addplot [domain=2:25,samples=1000,dashed]({1/x-1/x^2},{x});
    \end{axis}
\end{tikzpicture}%
\caption{ The backward Euler method for the equation $\dot q=q^2$ has
  a fold point at $(h_c,(2h_c)^{-1})$, where $h_c=(4x)^{-1}$. There
  are two solutions for $h<h_c$, a unique solution at $h=h_c$, and no
  solutions for $h>h_c$. Solutions on the lower branch, which we refer
  to as the \emph{principal solution branch}, are preferred because
  solutions they are naturally considered to be consistent
  with the initial condition.  }
\label{fig:demon-example}
\end{figure}

A simple example demonstrating the existence of a critical timestep $h_c$ at
which a fold bifurcation occurs is the application of the backward Euler
method to the scalar IVP
\begin{equation*}
  \dot{q} = q^2,\quad q(0) = q_0 > 0.
\end{equation*}
With the backward Euler method, the update
equation~\eqref{eq:general-time-integration1} can be written as
\begin{equation*}
  hy^2-y+x = 0,
\end{equation*}
having solutions
\begin{equation*}
  y = \frac{1\pm\sqrt{1-4hx}}{2h}.
\end{equation*}
Evidently, there are two solutions for $h<h_c=(4x)^{-1}$, a single solution at
$h=h_c$, and no solutions for $h>h_c$; see Figure~\ref{fig:demon-example}. We
note the existence of two solutions for all $0<h<h_c$. By the
Newton--Kantarovich Theorem~\cite{SchatzmanM-2002-1}, convergence to either
solution is possible for an appropriately chosen initial guess.

The principal solution branch contains the initial condition at zero
timestep, which is exactly correct, and it contains all solutions near
zero timestep that continuously emanate from the initial condition.
If a more complicated bifurcation occurs, say a pitchfork, as opposed
to a fold (e.g., bifurcations of the solutions of $y^3-(h-h_c)y=0$ and
$y^2-(h-h_c)=0$, respectively), then there are multiple solutions
beyond the critical timestep, none of which can be continuously and
monotonically continued back to the initial condition without passage
through the bifurcation itself. There may be two successive fold
bifurcations, resulting in an ``S''-shaped solution manifold, with
timesteps larger than $h_c$, for which there is a unique solution. But
in that situation, the two additional solutions at timesteps
$h\lesssim h_c$ would be rejected in favour of that on the principal
branch, and for $h\gtrsim h_c$, the unique solutions are near an
inconsistent one and so would also be rejected.

Therefore, if there are multiple solutions at a given timestep smaller
than the critical one, then the principal branch should be chosen, and
the timestep cannot be validly increased to be larger than the
critical one, irrespective of whether or not there are unique
solutions there.  Simulations should generally not be carried on from
solutions that are not on the principal branch, regardless of whether
the timestep used is larger or smaller the critical timestep, because
they tend to result in non-negligible perturbations to the solution. So
as to have a concise term, we make the following definitions.
\begin{definition}
  Given an initial condition $\bm x$, a solution $\bm y$ with
  timestep $h$ is called \emph{consistent} if $(\bm y, h)$ is on the principal
  solution branch.
\end{definition}

\begin{definition}
  The smallest timestep $h$ such that there is a bifurcation of the
  principal branch is called the critical timestep $h_c$. Timesteps
  $h > h_c$ are called \emph{invalid}. By definition, any solutions
  obtained with invalid timesteps cannot be on the principal branch.
\end{definition}

Another basic example relevant to Lagrangian systems that shows the
existence of a fold bifurcation for a critical timestep $h_c$ involves
the Lagrangian
\begin{equation*}
  \label{eq:lagrangian-example}
  L(x,\dot{x}) = \frac{1}{2}\dot{x}^2-\frac{1}{3}x^3,
\end{equation*}
which models the dynamics of a nonlinear spring that has a stiffness
proportional to its linear displacement from equilibrium. Applying the
first-order Variational Taylor method (described in
Section~\ref{sec:double-pendulum} below) yields an update equation of the form
\begin{equation*}
  \label{eq:example-VT}
  \begin{bmatrix}
    h\left(\dot{x}+\dot{y}\right) - 2\left(x-y\right) \\
    h\left(x^2+y^2\right) + 2\left(\dot{x}-\dot{y}\right)
  \end{bmatrix}
  = \mathbf{0},
\end{equation*}
to be solved for $y$ and $\dot{y}$.
The solution to this update equation is
\begin{equation}
  \label{eq:example-VT-update}
  \begin{aligned}
    y & =  h^{-2}\left(-2 \pm \sqrt{4-h^2\left(h^2x^2+4h\dot{x}-4x\right)}\right), \\
    \dot{y} & = -\dot{x} + 2h^{-1}\left(x-y\right),
  \end{aligned}
\end{equation}
which again shows the potential for multiple solutions. Taking, for example,
the initial state $x=1$, $\dot{x}=0$, it is straightforward to show that
(\ref{eq:example-VT-update}) has two real solutions for
$h<h_c=\sqrt{2+2\sqrt{2}}$, a single solution at $h=h_c$, and no real solutions
for $h>h_c$.


%% file: 0section-2-continuation.tex
\section{Numerical continuation of the implicit solution}
\label{sec:numerical-continuation}
We consider one-parameter families of solutions
to~\eqref{eq:general-time-integration1} in the hyperplane defined by
the vector $({\bm y},h)$ for fixed ${\bm x}$. For practical
computation, we parameterize these families by arclength, which
monotonically increases throughout the computation.

Suppose we have two points on the solution curve $({\bm y}_0,h_0)$ and
$({\bm y}_1,h_1)$ with known tangent vector ${\bm{\mathcal T}}_0$ at the first
point. Our goal is to find a next point on the solution curve in the same
arclength direction following $({\bm y}_1,h_1)$. This \emph{pseudo-arclength
  continuation} is accomplished in a predictor-corrector fashion. The
prediction is performed via a linear approximation to the curve at the point
$({\bm y}_1,h_1)$ and the correction by computing the minimum-norm solution to
the system with rank-1 deficiency.

The tangent vector $\bm{\mathcal T}_1$ at $({\bm y}_1,h_1)$ satisfies
\begin{equation*}
  \label{eq:tangent-1}
  \begin{bmatrix}
    \bm F_{\bm y} & \bm F_h
  \end{bmatrix}
  {\bm{\mathcal T}}_1 = \bm 0
\end{equation*}
and determines $\bm {\mathcal T}_1$ up to a scalar factor. To preserve the
direction of the orientation of the branch, we require
\begin{equation*}
  \label{eq:t-orientation}
  \bm {\mathcal T}_0^T\bm {\mathcal T}_1 = 1.
\end{equation*}
Using a natural decomposition
$\bm {\mathcal T}_i = (\bm {\mathcal T}_i^{(\bm y)},{\mathcal T}_i^{(h)})$, $i=0,1$, we write the
above as a single system
\begin{equation*}
  \label{eq:tangent}
  \begin{bmatrix}
    \bm F_{\bm y} & \bm F_h \\
    \bm {\mathcal T}_0^{(\bm y)T} & {\mathcal T}_0^{(h)}
  \end{bmatrix}
  \begin{bmatrix}
  \bm {\mathcal T}_1^{(\bm y)} \\
  {\mathcal T}_1^{(h)}
  \end{bmatrix}
  =
  \begin{bmatrix}
    \bm 0 \\
    1
  \end{bmatrix},
\end{equation*}
the solution to which uniquely determines $\bm {\mathcal T}_1$.

Using a step of length $\Delta s$, we create a predictor
\begin{equation*}
  \label{eq:predictor}
  \hat{\bm y}_2 = \bm y_1 + \frac{\Delta s}{\|\bm {\mathcal T}_1\|}\bm {\mathcal T}_1^{(\bm y)},\quad
  \hat{h}_2 = h_1 + \frac{\Delta s}{\|\bm {\mathcal T}_1\|}{\mathcal T}_1^{(h)},\quad
\end{equation*}
to approximate the next point on the curve. The name
\emph{pseudo-arclength} comes from the fact that $\Delta s$ measures
arclength along the tangent line.

From the predictor, a Newton-like method is applied to obtain the next
solution point on the curve $(\bm y_2,h_2)$. Our specific choice
of algorithm is the Gauss--Newton method. This particular approach for
path-following was first used in~\cite{Haselgrove1961}. It is
identified as being the Gauss--Newton method in~\cite{Deuflhard1987}.
A practical description is given in~\cite{Govaerts2000} and summarized
as follows.

We seek a solution to
\begin{equation*}
  \label{eq:Gauss-Newton-formulation}
  \bm F(\hat{\bm y}+\Delta\bm y;\bm x,\hat h+\Delta h) = \bm 0
\end{equation*}
such that $\|(\Delta\bm y,\Delta h)\|$ is minimal. Because $\bm F$
is nonlinear, we set up the Newton iteration as
$(\bm y^0,h^0)=(\hat{\bm y},\hat{h})$, and
$\bm y^{k+1} = \bm y^k + \Delta\bm y$, $h^{k+1} = h^k+\Delta h$, where
($\Delta\bm y$, $\Delta h$) is the minimum-norm solution to
\begin{equation*}
  \label{eq:Gauss-Newton-linear}
  \bm F_{\bm y}(\bm y^{k};\bm x,h^k)\Delta\bm y + \bm F_h(\bm y^{k};\bm x,h^k)\Delta h = -\bm F(\bm y^{k};\bm x,h^k).
\end{equation*}
The minimum-norm solution is obtained by solving the matrix system
\begin{equation*}
   \label{eq:Gauss-Newton-system}
  \begin{bmatrix}
    \bm F_{\bm y} & \bm F_h \\
    \bm {\mathcal T}_1^{(\bm y)T} & {\mathcal T}_1^{(h)}
  \end{bmatrix}
  \begin{bmatrix}
  \bm {\mathcal T}^{(\bm y)} &  \Delta_1\bm y \\
  {\mathcal T}^{(h)} &   \Delta_1 h

  \end{bmatrix}
   =
  \begin{bmatrix}
    \bm 0  &   -\bm F \\
    1 &    \hspace*{1.5ex} 0
  \end{bmatrix} \\
\end{equation*}
and constructing
\begin{equation*}
  \Delta\bm y = \Delta_1\bm y + \eta\bm {\mathcal T}^{(\bm y)},\quad
  \Delta h = \Delta_1 h + \eta {\mathcal T}^{(h)},
\end{equation*}
with
\begin{equation*}
  \label{eq:Gauss-Newton-eta}
  \eta = -\frac{(\Delta_1\bm y)^T\bm {\mathcal T}^{(\bm y)}+(\Delta_1
    h){\mathcal T}^{(h)}}{\|\bm {\mathcal T}\|^2}.
\end{equation*}
Quadratic convergence of the above algorithm is proven in
\cite{Deuflhard1987}, and in our experiments, we iterate solutions to a
tolerance $\|\bm F\| < F_{\text{tol}}=10^{-9}$, where
$\| \cdot \| = \| \cdot \|_2$.

An initial step length $\Delta s_0$ must be chosen to initialize the continuation.
For subsequent steps along the continuation curve, we choose the next step length
according to
\begin{equation*}
  \label{eq:steplength-control}
  \Delta s_{k} = \sqrt{\frac{2\epsilon}{\|\bm w_k\|}},\quad
    \bm w_k = \frac{1}{\Delta s_{k-1}}\left(\bm {\mathcal T}_k - \bm
      {\mathcal T}_{k-1}\right), \quad k=1,2,\ldots,
\end{equation*}
where $\epsilon$ is a user-defined tolerance for the absolute error in
$(\hat{\bm y},\hat h)$; we set $\epsilon=100\,F_{\text{tol}}$.

%% file: 0section-3-double-pendulum.tex
\section{The double pendulum}\label{sec:double-pendulum}
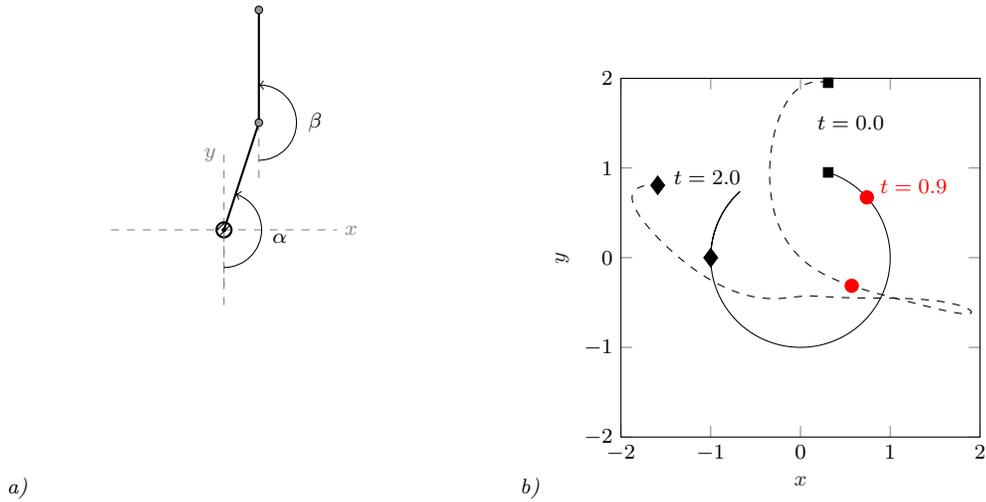
\begin{figure}[!htbp]
  \centering
  \begin{minipage}[t]{0.49\textwidth}
    \centering
    \footnotesize
\emph{a)}
  \begin{tikzpicture}[scale=0.5]

    \pgfmathsetmacro{\vert}{1.5}
    \pgfmathsetmacro{\myAnglea}{162}
    \pgfmathsetmacro{\myAngleb}{180}

    \draw[dashed,gray,-] (-3,0) -- (3,0) node (xaxis) [gray,right]{$x$};
    \draw[dashed,gray,-] (0,-2) -- (0,2) node (yaxis) [gray,left]{$y$};

    \coordinate (centro) at (0,0);
    \fill[pattern=north east lines] (centro) circle[radius=0.2];
    \draw[thick] (centro) circle[radius=0.2];
    \draw[dashed,gray,-] (centro) -- ++ (0,-\vert) node (mary) [black,below]{$ $};
    \draw[thick] (centro) -- ++(270+\myAnglea:3) coordinate (bob);
    \pic [draw, ->, "$\alpha$", angle eccentricity=1.5] {angle = mary--centro--bob};

    \draw[dashed,gray,-] (bob) -- ++ (0,-\vert) node (mary2) [black,below]{$ $};
    \draw[thick] (bob) -- ++(270+\myAngleb:3) coordinate (bob2);
    \pic [draw, ->, "$\beta$", angle eccentricity=1.5] {angle = mary2--bob--bob2};

    \filldraw [fill=black!40,draw=black] (bob) circle[radius=0.1];
    \filldraw [fill=black!40,draw=black] (bob2) circle[radius=0.1];
    \filldraw [fill=black,draw=black] (centro) circle[radius=0.05];

    \fill[fill=white] (0,-3) circle[radius=0.1];

    \clip (-5,-7) rectangle (5,5);
  \end{tikzpicture}
  \end{minipage}
  \begin{minipage}[t]{0.49\textwidth}
    \centering
    \setlength\fwidth{0.9\textwidth}
\footnotesize
\emph{b)} \input{figures/dp_trajectory_0_to_2.tex}
  \end{minipage}
  \caption{\emph{a)} Double pendulum in its initial configuration.
    \emph{b)} Reference trajectory of the double pendulum bobs used in
    this study. Starting from the initial configuration (black
    squares), the outer pendulum swings down performing a fast loop
    before arriving at the final configuration (black diamonds). We
    look at the behaviour of integrators starting from the
    configuration corresponding to time $t=0.9$ (red circles), just
    before entering the fast loop.}
  \label{fig:dp-viz}
\end{figure}

We applied pseudo-arclength continuation as described
in Section~\ref{sec:numerical-continuation} to the double pendulum system of two bobs
connected by a frictionless pin joint, the first moving on a circle with fixed
center, the second moving on a circle centered at the first, and both moving in
the presence of gravity; see Figure~\ref{fig:dp-viz}. The system is Lagrangian
(and so symplectic), with Lagrangian function
\begin{equation*}
  \label{eq:dp-polar-lagrangian}
  L^{\text{DP}}_{\text{ang}} = \dot\alpha^2+\frac{1}{2}\dot\beta^2+\dot\alpha\dot\beta\cos(\alpha-\beta) +
     g\left(2\cos\alpha+\cos\beta\right),
\end{equation*}
in terms of the angular coordinates $(\alpha,\beta)$ shown in
Figure~\ref{fig:dp-viz}a, and $g=9.81$.  The system evolution is in accord with
Euler--Lagrange differential equations of the form
\begin{equation*}
  \label{eq:EL-1storder}
  \dot{\bm q}(t)=\bm f(\bm q),\quad
   \bm q = \left(\alpha,\beta,\dot\alpha,\dot\beta\right)^T.
\end{equation*}
For initial conditions, we use the near-vertical arrangement as in Dharmaraja
et al.~\cite{Dharmaraja2010},
\begin{equation*} \label{eq:dp-initial-conditions}
    \alpha(0)=\frac{9\pi}{10}, \qquad
    \dot\alpha(0)=0.7,         \qquad
    \beta(0)=\pi,              \qquad
    \dot\beta(0)=0.4.
\end{equation*}

We compute a reference solution trajectory, $\bm q_{\text{ref}}$; details are
given at the end of this section. In the resulting motion, the system quickly
goes through a fast loop, as seen in the lower right of
Figure~\ref{fig:dp-viz}b. Events such as this, where the dynamics change
quickly relative to the initial or overall dynamics, seem to readily induce the
bifurcations in (\ref{eq:general-time-integration1}) in which we are
interested.  Accordingly, we chose a point on the trajectory that is just
before the loop for our $\bm x$ value, specifically the point corresponding to
time $t=0.9$.

The fold points we seek are a property of the implicit method used for
the integration as defined by $\bm F(\bm y;\bm x, h)$.  We now
investigate the properties of these fold points for various implicit
integration methods.

We first consider a Variational Taylor~(\VarTay)
method~\cite{Marsden2001,Patrick2009}, which uses the first-order Taylor
expansion of trajectories $ \bm x \rightarrow \bm x +t\dot{\bm x} $ to obtain
the \emph{discrete Lagrangian},
\[
L_h(\bm x) = \int_{-h/2}^{h/2}L^{DP}_{\text{ang}}(\bm x+t\dot{\bm x})\,\mathrm{d}t.
\]
The resulting one-step method is obtained by finding the
critical points of the \emph{discrete action} $L_h(\bm y)+L_h(\bm x)$
subject to the constraints
\[
  \begin{bmatrix}
    x_1 \\ x_2
  \end{bmatrix}
-\frac{h}{2}
\begin{bmatrix}
  \dot{x}_1 \\ \dot{x}_2
\end{bmatrix}
= \bm a_1, \quad
  \begin{bmatrix}
    y_1 \\ y_2
  \end{bmatrix}
+\frac{h}{2}
\begin{bmatrix}
  \dot{y}_1 \\ \dot{y}_2
\end{bmatrix}
= \bm a_2, \quad
  \begin{bmatrix}
    x_1 \\ x_2
  \end{bmatrix}
+\frac{h}{2}
\begin{bmatrix}
  \dot{x}_1 \\ \dot{x}_2
\end{bmatrix}
=
  \begin{bmatrix}
    y_1 \\ y_2
  \end{bmatrix}
-\frac{h}{2}
\begin{bmatrix}
  \dot{y}_1 \\ \dot{y}_2
\end{bmatrix},
\]
where $\bm a_1$ and $\bm a_2$ are constants, the values of which are
not required for the method implementation. The first two of these
constraints are discrete analogues of the fixed-endpoint constraint
for the continuous variational principle. The last, which connects the
discrete trajectories, is a discrete analogue of the continuous
constraint that the derivative of configuration is velocity. The
resulting constrained optimization problem is equivalent to
\emph{discrete Euler--Lagrange equations}, and, in the case at hand,
the associated Lagrange multipliers may be explicitly eliminated to
obtain equations of the
form~\eqref{eq:general-time-integration1}. Such methods extend to any
order, but they become much more complicated, and it suffices for our
purpose here to consider only the first-order
method~\cite{PatrickGW-CuellC-SpiteriRJ-ZhangW-2009-1,PatrickGW-2009-1}.

In addition to this, we investigate the fold points of a variety of
implicit Runge--Kutta methods of the form
\begin{equation*}
\label{eq:RK-method}
  \bm q^{n+1} = \bm q^n + h\sum_{i=1}^{s}b_i \bm k_i,\qquad
  \bm k_i = \bm f\left(t_n+c_ih,\bm q_n+h\sum_{j=1}^sa_{ij}\bm k_j\right).
\end{equation*}
Specifically, using the usual Butcher tableau
notation~\cite{HairerNorsettWanner1993}, we use the methods below (in
row order): backward Euler (\BE); trapezoidal rule (\CN); the
trapezoidal rule backward differentiation formula 2 split-step method,
with optimized split-step size $\gamma=2-\sqrt{2}$ as in
\cite{Dharmaraja2010} (\trbdf); third- and fifth-order Radau IIA
methods (\RadThree, \RadFive); and fourth- and sixth-order
Gauss--Legendre methods (\GauLegFour, \GauLegSix).  In the experiments
described below, the nonlinear systems associated with any implicit
method were solved using a classical Newton iteration.

\begin{center}
\parbox{0cm}{
\footnotesize\renewcommand{\arraystretch}{1.75}
\begin{tabbing}
AAAAAAAA\=AAAAAAAAAAAAA\=\kill
$\input{0tableau-backward-euler.tex}$
\>
$\input{0tableau-crank-nicolson.tex}$
\>
$\input{0tableau-trbdf2.tex}$
\\\\
$\input{0tableau-radauIIA-order-3.tex}$
\>\>
$\input{0tableau-radauIIA-order-5.tex}$
\\\\
$\input{0tableau-gauss-legendre-order-4.tex}$
\>\>
$\input{0tableau-gauss-legendre-order-6.tex}$
\end{tabbing}}\end{center}

A reference solution for this problem was computed using timestep
$h_{\text{ref}}=2\times 10^{-5}$ with the \VarTay\ method. This
reference solution was compared to that obtained from a \VarTay\
integration with timestep size $h_{\text{ref}}/2$. The state of the
system at $t=2.0$ is given in Table~\ref{tab:ref-digits}, where we see
agreement of 8 significant digits between these solutions.

\begin{table}[!htbp]
  \centering
  \begin{tabular}{r|r|r|c}
    & \multicolumn{1}{c|}{$h_{\text{ref}}$} &
        \multicolumn{1}{c|}{$h_{\text{ref}}/2$} & \bf Digits \\
    \hline
    $\alpha$ &$ \bm{-1.5707374}51319846$ & $\bm{-1.5707374}39361913$ & $8$\\
    $\beta$   & $\bm{3.7730189}50076258$  & $\bm{3.7730189}44217077$  & $8$\\
    $\dot{\alpha}$ & $\bm{4.1181166}60671203$  & $\bm{4.1181166}38219623$  & $8$\\
    $\dot{\beta}$   & $\bm{-6.2736260}26547350$ & $\bm{-6.2736260}00367146$ & $8$
  \end{tabular}
  \caption{State of the reference solution at final time $t=2.0$ compared
    to that obtained by halving the timestep. Matching digits are in
    bold font.}
  \label{tab:ref-digits}
\end{table}

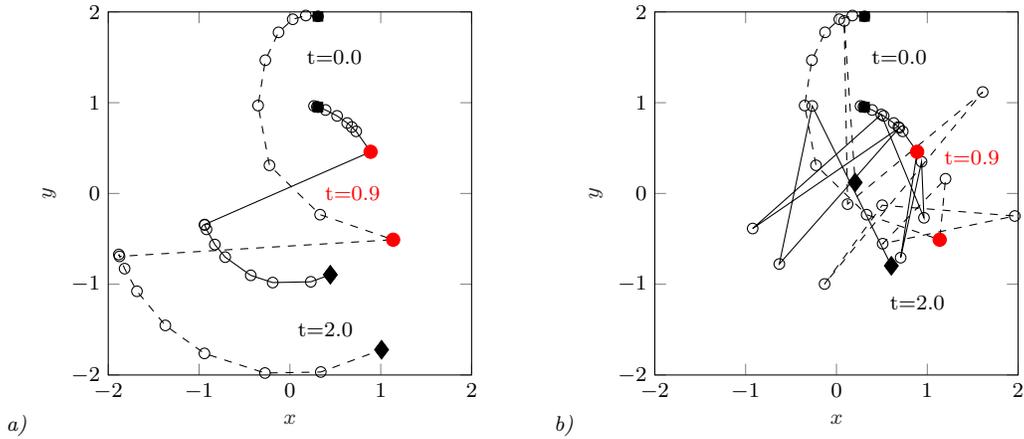
\begin{figure}[!htbp]
  \centering
  \begin{minipage}[t]{0.49\textwidth}
    \centering
    \setlength\fwidth{0.9\textwidth}
\footnotesize
\emph{a)} \input{figures/dp_bad_trajectory_noextrap.tex}
  \end{minipage}
  \begin{minipage}[t]{0.49\textwidth}
    \centering
    \setlength\fwidth{0.9\textwidth}
\footnotesize
\emph{b)} \input{figures/dp_bad_trajectory_extrap.tex}
  \end{minipage}
  \caption{Simulations of the double pendulum using a timestep of
    $h=0.1225$ with the VT method. \emph{a)} The nonlinear solver of a
    VT timestep (initialized with the solution at the previous timestep)
    converges to a solution, but the trajectory differs
    substantially from the true trajectory as shown in Figure
    \ref{fig:dp-viz}. \emph{b)} Changing the initialization of the
    nonlinear solver to a linear extrapolation based on the solution at
    the previous two timesteps, the simulation still converges, but
    the resulting trajectory diverges even further.}
  \label{fig:dp-bad-trajectories}
\end{figure}

To demonstrate the potential effects of convergence to
non-principal-branch solutions from an implicit one-step method, we
show inconsistent solutions obtained with the VT method and timestep
size $h=0.1225$ in Figure~\ref{fig:dp-bad-trajectories}. The nonlinear
solver used in computing the solution at the next timestep converges
in all steps displayed, yet the inconsistent solutions appear markedly
different from the reference solution of
Figure~\ref{fig:dp-viz}. There is even a difference in trajectories
with inconsistent solutions when all that is changed is the
initialization of the nonlinear solver, in this case, from the
solution at the previous timestep to a linear extrapolation based on
the solution at the previous two timesteps; this implies the existence
of at least two non-principal solution branches at some point along
the time integration.

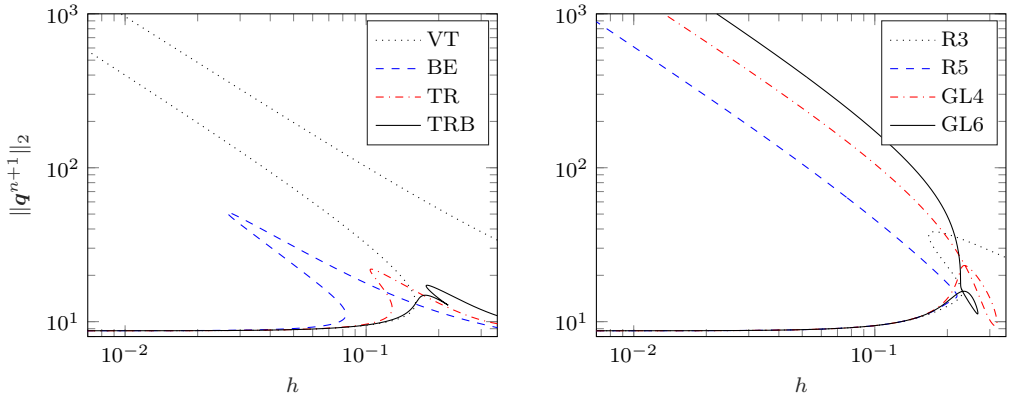
\begin{figure}[htbp!]
  \centering
  \begin{minipage}{0.49\textwidth}
    \setlength\fwidth{0.8\textwidth}
    \footnotesize
    \input{figures/low_order_norms.tex}
  \end{minipage}
  \begin{minipage}{0.49\textwidth}
    \setlength\fwidth{0.8\textwidth}
    \footnotesize
    \input{figures/high_order_norms.tex}
  \end{minipage}
  \caption{Norm of the computed state for various integrators. Point
    on the trajectory used for initialization is at $t=0.9$ from the
    reference solution pictured in Figure \ref{fig:dp-viz}b. }
  \label{fig:dp-norm-state}
\end{figure}

All of the methods tested exhibit fold points before $h=0.33$, as seen
in Figure~\ref{fig:dp-norm-state}.  Four of the methods (\BE, \CN,
\trbdf, and \RadThree) have two fold points, resulting in a continuous
connection between the $h=0$ solution and the $h=0.35$ solution.
Three of the methods (\RadFive, \GauLegFour, and \GauLegSix) have a
single fold point, beyond which we were unable to find other
solutions.  The solution obtained with the \VarTay\ method exhibits
two separate solution branches.  The main branch originating from the
$h=0$ solution folds back, with norm approaching infinity as $h$
returns to zero.  A second branch extends beyond the fold point,
allowing for solutions for larger $h$ that are not continuously and monotonically
connected to the $h=0$ solution.

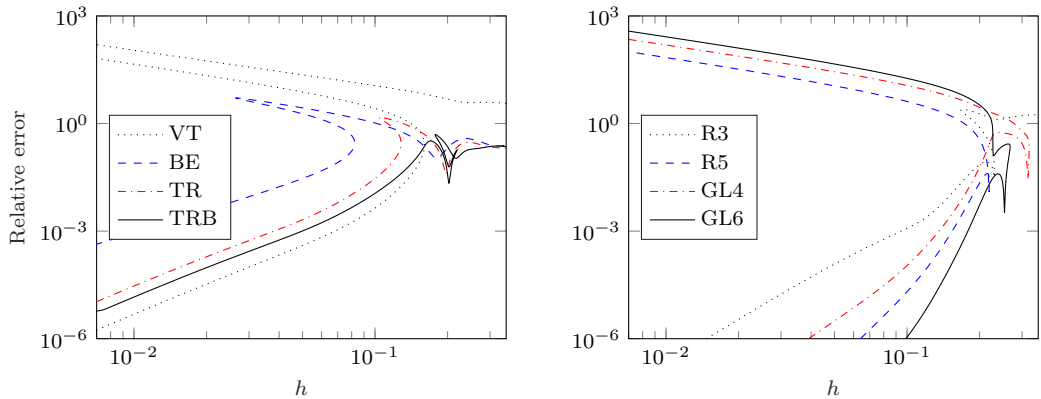
\begin{figure}[htbp!]
  \centering
  \begin{minipage}{0.49\linewidth}
    \centering
    \setlength\fwidth{0.8\textwidth}
    \footnotesize
    \input{figures/low_order_traj_errors.tex}
  \end{minipage}
  \begin{minipage}{0.49\linewidth}
    \centering
    \setlength\fwidth{0.8\textwidth}
    \footnotesize
    \input{figures/high_order_traj_errors.tex}
  \end{minipage}
  \caption{Relative errors in the state $\bm q^{n+1}$.  These errors
    are large (typically $\mathcal{O}(1)$) at the fold points.  }
  \label{fig:dp-errors}
\end{figure}

The relative error in the trajectories, computed as
\[
e_{rel} = \frac{\|\bm q^{n+1}-\bm q_{\text{ref}}^{n+1}\|_2}{\|\bm q_{\text{ref}}^{n+1}\|_2},
\]
is visualized for all solution curves in Figure \ref{fig:dp-errors}.
We notice that the computed solutions near the fold points for each of
the methods have relatively large (typically $\mathcal{O}(1)$) errors
in this instance and that the methods are out of their regions of
asymptotic convergence.  Passing through the fold points in arclength,
the error generally increases as the timestep decreases.

%% file: figures/dp_trajectory_0_to_2.tex
%
%
\begin{tikzpicture}

\begin{axis}[%
width=0.741\fwidth,
height=0.741\fwidth,
at={(0\fwidth,0\fwidth)},
scale only axis,
xmin=-2,
xmax=2,
xlabel style={font=\color{white!15!black}},
xlabel={$x$},
ymin=-2,
ymax=2,
ylabel style={font=\color{white!15!black}},
ylabel={$y$},
axis background/.style={fill=white}
]
\addplot [color=black, forget plot]
  table[row sep=crcr]{%
0.309016994374948	0.951056516295154\\
0.263490849581812	0.964661895270387\\
0.309503788676464	0.950898209481391\\
0.354807870200931	0.934939236123653\\
0.399287674786106	0.916825693773852\\
0.442846147449762	0.896597618605416\\
0.485394841012255	0.87429505792878\\
0.526819204005815	0.849977368104986\\
0.567024782433796	0.823700732126606\\
0.605927164899611	0.795520125978419\\
0.643433135329112	0.765502318977904\\
0.679451874306461	0.733720076392514\\
0.691539006295932	0.722339118953996\\
0.690451472084978	0.723378714571895\\
0.724387032615598	0.689393521131704\\
0.756642639045591	0.653828660107618\\
0.787171067695275	0.616734716213933\\
0.81586924739452	0.578236431882236\\
0.842698912632848	0.538385124838545\\
0.867552238945224	0.497346069353251\\
0.89042458262174	0.455130819284852\\
0.911257598741292	0.411836847229888\\
0.929979478174171	0.367611439124106\\
0.946558242074079	0.322532935313013\\
0.960944664866894	0.276740584417706\\
0.973118259684726	0.230305997898817\\
0.983048945830442	0.183343312126887\\
0.990717712165499	0.135935333160877\\
0.996097271763903	0.0882622523194865\\
0.999183973208042	0.0403904405050386\\
0.999970999960916	-0.00761572302318525\\
0.998448778961101	-0.0556779650408205\\
0.994625977717788	-0.103533397746499\\
0.988509613188295	-0.151158012140698\\
0.980103032201358	-0.198489410976262\\
0.969459392861698	-0.245251881934122\\
0.95658441259997	-0.291455419525474\\
0.941527659542137	-0.336935700567787\\
0.924283472001372	-0.381706776708365\\
0.904884007749585	-0.425658234407662\\
0.883413623544267	-0.468594035105429\\
0.859931404148205	-0.51040961997174\\
0.834453112553609	-0.551078944389635\\
0.807031375429121	-0.59050855969493\\
0.777767678069503	-0.628551858600683\\
0.746730425702642	-0.66512680845832\\
0.713955506028932	-0.700191070644986\\
0.679545719519619	-0.733633161111576\\
0.643585161612784	-0.76537450947353\\
0.606142602179756	-0.795355986852903\\
0.56727516206966	-0.823528317970209\\
0.527031979319467	-0.849845452288005\\
0.485670422250946	-0.874142002737878\\
0.443139191747975	-0.896452819024488\\
0.399603459346512	-0.916688101416343\\
0.355131871855547	-0.934816213804818\\
0.3098042806434	-0.950800351122687\\
0.263833386018636	-0.96456826840921\\
0.217218271979549	-0.976123056954511\\
0.170091602504315	-0.985428255510016\\
0.122608453799496	-0.992455120928345\\
0.0748191905490737	-0.997197116284229\\
0.0269192691061729	-0.999637610812433\\
-0.0211377658880827	-0.999776572466699\\
-0.0691229615936192	-0.997608147611339\\
-0.116922662136692	-0.993141022755011\\
-0.164411715670135	-0.98639180235361\\
-0.211593904582399	-0.977357672269254\\
-0.258190323525445	-0.96609407245766\\
-0.304193035870154	-0.952610411935592\\
-0.349591524292641	-0.936902218026379\\
-0.394095642076635	-0.919069434208431\\
-0.437692846500383	-0.899124558736103\\
-0.48036849507856	-0.877066764240875\\
-0.521971509392776	-0.852962920285652\\
-0.562355645530138	-0.82689547582532\\
-0.601381995703373	-0.798961635652069\\
-0.639043139481659	-0.769170895108119\\
-0.675208222272025	-0.737627179933232\\
-0.709869613055242	-0.704333111858872\\
-0.742903327953068	-0.66939871923709\\
-0.774196387662915	-0.6329454584162\\
-0.803743919707573	-0.594975387333885\\
-0.831350177613239	-0.555748938084847\\
-0.857106817659208	-0.515138722211896\\
-0.880837574063199	-0.473418597140479\\
-0.902544958588843	-0.430595631336249\\
-0.922155633243481	-0.386819063745978\\
-0.939664768893842	-0.342096656077901\\
-0.95500606145251	-0.296586281862403\\
-0.968128660246018	-0.250453383307652\\
-0.979029452590861	-0.203718754560396\\
-0.987668481413276	-0.15655979951058\\
-0.994040540955734	-0.109011022086902\\
-0.99813054769029	-0.0611179987195453\\
-0.999914359439862	-0.0130871611119722\\
-0.999391956614913	0.0348671342874121\\
-0.996565814317076	0.0828044547989035\\
-0.991433211978652	0.13061464766899\\
-0.984024691800785	0.178032036236093\\
-0.97435004239354	0.225037763248096\\
-0.962419464749674	0.271567254784059\\
-0.948294829005813	0.317390795835095\\
-0.931954763623111	0.362574569654563\\
-0.913484135417395	0.406874347115585\\
-0.892893363105761	0.450268189106985\\
-0.870263602842497	0.492586298599135\\
-0.84560360251978	0.533811340649081\\
-0.819016987749169	0.573769268764263\\
-0.790525144763386	0.612429584113004\\
-0.760203284622683	0.649685282308968\\
-0.72810804596552	0.685462379274218\\
-0.694356660000343	0.719631036513274\\
-0.674437203351162	0.738332214342476\\
-0.709098630404831	0.70510930525557\\
-0.742135504936971	0.670249873041351\\
-0.773473012258719	0.633829235131532\\
-0.802991149879916	0.59599095061463\\
-0.830666861047543	0.556769760275666\\
-0.856415718150311	0.516286856025879\\
-0.880231255672906	0.474544978412478\\
-0.902001760935167	0.431732351428356\\
-0.921677225601575	0.387957590217518\\
-0.93923024969866	0.343287835570959\\
-0.954608907612557	0.297862104851827\\
-0.967808898600785	0.251686185137634\\
-0.97878094066303	0.204909419487721\\
-0.98748591724697	0.157707207314412\\
-0.993911845993954	0.110178230122337\\
-0.998053456293855	0.062364239592015\\
-0.999896890445332	0.0143599609245879\\
-0.999999998266839	-5.88754750170395e-05\\
};
\addplot [color=black, dashed, forget plot]
  table[row sep=crcr]{%
0.309016994374948	1.95105651629515\\
0.274233185371354	1.95692866284925\\
0.23921093071993	1.96049010456345\\
0.205104392279125	1.96163837610802\\
0.172614416318772	1.96042512292138\\
0.141903812469519	1.95696471455331\\
0.112842732357765	1.95135420393125\\
0.0851575878502484	1.94363359999678\\
0.0585909773826945	1.93379388629156\\
0.0328905446110348	1.9217695562914\\
0.00787524066923484	1.90746803314397\\
-0.0166127066938124	1.89075656048207\\
-0.0407076568241755	1.87146020282151\\
-0.0644683366297276	1.8494114802119\\
-0.0879246503193118	1.82443132323043\\
-0.111075988150593	1.79633751011578\\
-0.133864862945846	1.76499194794678\\
-0.156251151402316	1.73021330981887\\
-0.178112491117907	1.69193034196712\\
-0.199343791799223	1.65004576673661\\
-0.219808215261425	1.60451912508436\\
-0.239318295204259	1.55544419274953\\
-0.257707581429913	1.50290547165519\\
-0.274768797154616	1.44715479253983\\
-0.290321805486984	1.38844114329302\\
-0.304177376439672	1.32713777185779\\
-0.316176742979226	1.26360219904776\\
-0.326151195780017	1.19840421003137\\
-0.333986072886891	1.1320116400859\\
-0.339584333756038	1.06502117462641\\
-0.342896600762909	0.99788618434978\\
-0.343898403447856	0.931172226625056\\
-0.342601920709646	0.865263291223808\\
-0.339041673159067	0.800535771462365\\
-0.333278472590113	0.737355751513572\\
-0.325387480123218	0.675966349606791\\
-0.315439142850023	0.616474786868737\\
-0.303527273262178	0.559068722109126\\
-0.289700901125628	0.503684181624763\\
-0.274055555308884	0.450452933094682\\
-0.256659599371734	0.399369234040249\\
-0.237525522167936	0.350290257502197\\
-0.216737019648913	0.303280133021542\\
-0.194328028486477	0.258271932878765\\
-0.170316333565721	0.215184722605738\\
-0.144701544986704	0.173927346421352\\
-0.11746355117507	0.134402296623129\\
-0.0885615740770858	0.0965095215389145\\
-0.0579339208378815	0.0601500009878806\\
-0.0254984884402627	0.025228900830979\\
0.00874826565332532	-0.00825085349766241\\
0.0450125652336741	-0.0404682155847389\\
0.0833208807492756	-0.071411669656843\\
0.123917015929416	-0.10123483829804\\
0.167070700062199	-0.130066258637084\\
0.213070552948513	-0.158011268506834\\
0.262216741589088	-0.185156281267399\\
0.315064005285789	-0.211695305581742\\
0.372183974798811	-0.237784159889197\\
0.434539847484589	-0.263709757732179\\
0.503366560645204	-0.28980343740287\\
0.580556672987677	-0.316569989492818\\
0.669183335363045	-0.344818938542731\\
0.774668656030849	-0.375963717674282\\
0.909070538789302	-0.413146622648896\\
1.10682828896376	-0.465272849787138\\
1.40040561276803	-0.540275971941194\\
1.55763706578789	-0.578393084867356\\
1.65735912698657	-0.600354234264297\\
1.72950576761303	-0.614018586402362\\
1.78379412084291	-0.622074305534357\\
1.82466342465527	-0.625925806407132\\
1.85498134295017	-0.626573198419957\\
1.87645957962604	-0.624830864468548\\
1.89077761199571	-0.621445734699751\\
1.89942854733445	-0.617072072017437\\
1.90382202501419	-0.612186905218847\\
1.90507830712784	-0.606850388059138\\
1.90353302541244	-0.600802746225461\\
1.89871975563247	-0.593658990749606\\
1.88970784551653	-0.585195959570451\\
1.87545661561232	-0.575376808701346\\
1.85488930283731	-0.564263417809191\\
1.8274074071622	-0.552200455257789\\
1.79255129164653	-0.539529383736548\\
1.75016907327011	-0.526652827768509\\
1.70029467050069	-0.513963892473461\\
1.64296451022268	-0.501801292475049\\
1.57853799680653	-0.490517594213913\\
1.50692389521703	-0.480335689980261\\
1.42804727866076	-0.471465684499813\\
1.34132867830556	-0.464050420087232\\
1.24556916736593	-0.458196440369757\\
1.1380800699391	-0.453961379714057\\
1.01127206121892	-0.451327761770992\\
0.829316049713466	-0.450081422788921\\
0.59826508305423	-0.448137927596752\\
0.476026875937266	-0.444736901495412\\
0.363170657465043	-0.439296985317694\\
0.125787980956546	-0.426715504854091\\
0.0723824507047992	-0.427202557971402\\
0.021205217070823	-0.429995843430795\\
-0.0380987665859545	-0.435663169688884\\
-0.210113864618585	-0.453630295401128\\
-0.266686371128557	-0.456336115391479\\
-0.319990602212721	-0.456595278437122\\
-0.371973200657077	-0.454550727891917\\
-0.423488269117967	-0.450222474389725\\
-0.475110175724184	-0.443568849209933\\
-0.526996637732961	-0.434552001472436\\
-0.579262617235305	-0.423129663573613\\
-0.63219619654876	-0.40919690480233\\
-0.685765387248599	-0.392706923976211\\
-0.739955404737721	-0.373614162086858\\
-0.794887091582616	-0.351819779678829\\
-0.850688802201393	-0.327200415561145\\
-0.907366758118269	-0.299672169900417\\
-0.964923024518893	-0.269152532568139\\
-1.02335231219781	-0.235563101399668\\
-1.08275411295712	-0.198758806648844\\
-1.14309181138425	-0.158668702080354\\
-1.20441991472433	-0.115157709070469\\
-1.26654472236892	-0.0682650474243449\\
-1.32935910631675	-0.0179766893503366\\
-1.39251925811367	0.0355272630512591\\
-1.45527904706022	0.0916868011034777\\
-1.51682971290733	0.149822374787502\\
-1.57587409476345	0.208715189360416\\
-1.63091479752767	0.266805548461338\\
-1.68073877375673	0.322666306440131\\
-1.7243489875763	0.37493264573295\\
-1.76153301298628	0.423002104904188\\
-1.79233656850528	0.466487511789629\\
-1.81730922926181	0.505629219612293\\
-1.83699380056859	0.54067106089176\\
-1.85201328921692	0.572024509665152\\
-1.86294064665166	0.60009388720052\\
-1.8702902617405	0.625275198171767\\
-1.87450446065188	0.647952876355992\\
-1.87593673860213	0.668350624032972\\
-1.87489364662148	0.686829650244688\\
-1.87156687980804	0.703717584165721\\
-1.86607767208542	0.719176526839403\\
-1.85849064474582	0.733346679415443\\
-1.84872454967716	0.746440728033897\\
-1.83672948195807	0.758461932689888\\
-1.82225807739219	0.769513039576165\\
-1.80515643458552	0.77950770750899\\
-1.7851895556177	0.788363993211076\\
-1.76217139091181	0.795931045855656\\
-1.73591575349505	0.80203456499129\\
-1.70623763445731	0.80647676412341\\
-1.67314143484417	0.809027450127221\\
-1.63656654566351	0.80947630649164\\
-1.59656257671019	0.807615701956252\\
-1.59029664335955	0.807127516129973\\
};
\addplot [color=black, draw=none, mark size=1.8pt, mark=square*, mark options={solid, fill=black, black}, forget plot]
  table[row sep=crcr]{%
0.309016994374948	0.951056516295154\\
};
\addplot [color=black, draw=none, mark size=1.8pt, mark=square*, mark options={solid, fill=black, black}, forget plot]
  table[row sep=crcr]{%
0.309016994374948	1.95105651629515\\
};
\node[right, align=left]
at (axis cs:0.1,1.5) {$t=0.0$};
\addplot [color=red, draw=none, mark size=2.5pt, mark=*, mark options={solid, fill=red, red}, forget plot]
  table[row sep=crcr]{%
0.740357312510513	0.672213544799724\\
};
\addplot [color=red, draw=none, mark size=2.5pt, mark=*, mark options={solid, fill=red, red}, forget plot]
  table[row sep=crcr]{%
0.570773769739565	-0.31330227040633\\
};
\node[right, align=left, font=\color{red}]
at (axis cs:0.8,0.8) {$t=0.9$};
\addplot [color=black, draw=none, mark size=3.5pt, mark=diamond*, mark options={solid, fill=black, black}, forget plot]
  table[row sep=crcr]{%
-0.999999998266839	-5.88754750170395e-05\\
};
\addplot [color=black, draw=none, mark size=3.5pt, mark=diamond*, mark options={solid, fill=black, black}, forget plot]
  table[row sep=crcr]{%
-1.59029664335955	0.807127516129973\\
};
\node[right, align=left]
at (axis cs:-1.5,0.9) {$t=2.0$};
\end{axis}
\end{tikzpicture}%

%% file: 0tableau-backward-euler.tex
\begin{array}{c|c}
  1 & 1
  \\
  \hline
  & 1
\end{array} 

%% file: 0tableau-crank-nicolson.tex
\begin{array}{c|cc}
  0 & 0 & 0
  \\
  1 & 1/2 & 1/2
  \\
  \hline
  & 1/2 & 1/2
\end{array} 

%% file: 0tableau-trbdf2.tex
\begin{array}{c|ccc}
  0 & 0 & 0 & 0
  \\
  \gamma & \frac{\gamma}{2} & \frac{\gamma}{2} & 0
  \\
  1
  & \frac{1}{2(2-\gamma)}
  & \frac{1}{2(2-\gamma)}
  & \frac{1-\gamma}{2-\gamma}
  \\
  \hline
  & \frac{1}{2(2-\gamma)}
  & \frac{1}{2(2-\gamma)}
  & \frac{1-\gamma}{2-\gamma}
  \\
\end{array} 

%% file: 0tableau-radauIIA-order-3.tex
\begin{array}{c|cc}
  1/3 & 5/12 & -1/12
  \\
  1 & 3/4 & 1/4
  \\
 \hline
  & 3/4 & 1/4
\end{array} 

%% file: 0tableau-radauIIA-order-5.tex
\begin{array}{c|ccc}
  \frac{2}{5}-\frac{\sqrt{6}}{10}
  & \frac{11}{45}-\frac{7\sqrt{6}}{360}
  & \frac{37}{225}-\frac{169\sqrt{6}}{1800}
  & -\frac{2}{225}+\frac{\sqrt{6}}{75}
  \\
  \frac{2}{5} + \frac{\sqrt{6}}{10}
  & \frac{37}{225}+\frac{169\sqrt{6}}{1800}
  & \frac{11}{45}+\frac{7\sqrt{6}}{360}
  & -\frac{2}{225}-\frac{\sqrt{6}}{75}
  \\
  1
  & \frac{4}{9}-\frac{\sqrt{6}}{36}
  & \frac{4}{9}+\frac{\sqrt{6}}{36}
  & \frac{1}{9}
  \\
  \hline
  & \frac{4}{9}-\frac{\sqrt{6}}{36}
  & \frac{4}{9}+\frac{\sqrt{6}}{36}
  & \frac{1}{9}
\end{array} 

%% file: 0tableau-gauss-legendre-order-4.tex
\begin{array}{c|cc}
  \frac{1}{2}-\frac{\sqrt{3}}{6}
  & \frac{1}{4}
  & \frac{1}{4}-\frac{\sqrt{3}}{6}
  \\
  \frac{1}{2}+\frac{\sqrt{3}}{6}
  & \frac{1}{4}+\frac{\sqrt{3}}{6}
  & \frac{1}{4}
  \\
  \hline
  & \frac{1}{2}
  & \frac{1}{2}
\end{array}

%% file: 0tableau-gauss-legendre-order-6.tex
\begin{array}{c|ccc}
  \frac{1}{2}-\frac{\sqrt{15}}{10}
  & \frac{5}{36}
  & \frac{2}{9}-\frac{\sqrt{15}}{15}
  & \frac{5}{36}-\frac{\sqrt{15}}{30}
  \\
  \frac{1}{2}
  & \frac{5}{36}+\frac{\sqrt{15}}{24}
  & \frac{2}{9}
  & \frac{5}{36}-\frac{\sqrt{15}}{24}
  \\
  \frac{1}{2}+\frac{\sqrt{15}}{10}
  & \frac{5}{36}+\frac{\sqrt{15}}{30}
  & \frac{2}{9}+\frac{\sqrt{15}}{15}
  & \frac{5}{36}
  \\
  \hline
  & \frac{5}{18}
  & \frac{4}{9}
  & \frac{5}{18}
\end{array}

%% file: figures/dp_bad_trajectory_noextrap.tex
%
%
\begin{tikzpicture}

\begin{axis}[%
width=0.75\fwidth,
height=0.75\fwidth,
at={(0\fwidth,0\fwidth)},
scale only axis,
xmin=-2,
xmax=2,
xlabel style={font=\color{white!15!black}},
xlabel={$x$},
ymin=-2,
ymax=2,
ylabel style={font=\color{white!15!black}},
ylabel={$y$},
axis background/.style={fill=white}
]
\addplot [color=black, mark=o, mark options={solid, black}, forget plot]
  table[row sep=crcr]{%
0.309016994374948	0.951056516295154\\
0.265115830725854	0.964216571263189\\
0.296851532460923	0.954923644945292\\
0.392718721303604	0.919658635547811\\
0.519296828849877	0.854593940738209\\
0.630793525981773	0.77595072496872\\
0.681447508104726	0.731866991807842\\
0.729428687256807	0.684056861822766\\
0.888386102650158	0.459097084087955\\
-0.939391135477632	-0.34284733422625\\
-0.937636650160155	-0.347616904474513\\
-0.91840986876814	-0.395630273044527\\
-0.826543463124387	-0.562872901787201\\
-0.71241500388876	-0.701758407312786\\
-0.431263329687403	-0.902226102740845\\
-0.189805874345477	-0.981821638620757\\
0.229998449101707	-0.973190995339974\\
0.444860921775787	-0.8955996651835\\
};
\addplot [color=black, dashed, mark=o, mark options={solid, black}, forget plot]
  table[row sep=crcr]{%
0.309016994374948	1.95105651629515\\
0.175320459892549	1.96017680710441\\
0.0310067596273049	1.91893949375538\\
-0.125349647149275	1.77499784639259\\
-0.270533358804915	1.46791952692323\\
-0.350391159922619	0.969022245775647\\
-0.225191669853145	0.309960128872578\\
0.332330239226565	-0.233719158005587\\
1.13637255651907	-0.509666415790897\\
-1.87599574324481	-0.693235417227488\\
-1.88347408819203	-0.672257537817438\\
-1.81992965592485	-0.8283681537124\\
-1.68382391701864	-1.07772261128587\\
-1.37142087976961	-1.4538961975843\\
-0.941421326105866	-1.76230680749494\\
-0.274752154528746	-1.97820717117481\\
0.340144711058563	-1.96710638458442\\
1.00875421899317	-1.72144732023455\\
};
\addplot [color=black, draw=none, mark size=1.8pt, mark=square*, mark options={solid, fill=black, black}, forget plot]
  table[row sep=crcr]{%
0.309016994374948	0.951056516295154\\
};
\addplot [color=black, draw=none, mark size=1.8pt, mark=square*, mark options={solid, fill=black, black}, forget plot]
  table[row sep=crcr]{%
0.309016994374948	1.95105651629515\\
};
\node[right, align=left]
at (axis cs:0.1,1.5) {t=0.0};
\addplot [color=red, draw=none, mark size=2.5pt, mark=*, mark options={solid, fill=red, red}, forget plot]
  table[row sep=crcr]{%
0.888386102650158	0.459097084087955\\
};
\addplot [color=red, draw=none, mark size=2.5pt, mark=*, mark options={solid, fill=red, red}, forget plot]
  table[row sep=crcr]{%
1.13637255651907	-0.509666415790897\\
};
\node[right, align=left, font=\color{red}]
at (axis cs:0.3,0) {t=0.9};
\addplot [color=black, draw=none, mark size=3.5pt, mark=diamond*, mark options={solid, fill=black, black}, forget plot]
  table[row sep=crcr]{%
0.444860921775787	-0.8955996651835\\
};
\addplot [color=black, draw=none, mark size=3.5pt, mark=diamond*, mark options={solid, fill=black, black}, forget plot]
  table[row sep=crcr]{%
1.00875421899317	-1.72144732023455\\
};
\node[right, align=left]
at (axis cs:0,-1.5) {t=2.0};
\end{axis}
\end{tikzpicture}%

%% file: figures/dp_bad_trajectory_extrap.tex
%
%
\begin{tikzpicture}

\begin{axis}[%
width=0.75\fwidth,
height=0.75\fwidth,
at={(0\fwidth,0\fwidth)},
scale only axis,
xmin=-2,
xmax=2,
xlabel style={font=\color{white!15!black}},
xlabel={$x$},
ymin=-2,
ymax=2,
ylabel style={font=\color{white!15!black}},
ylabel={$y$},
axis background/.style={fill=white}
]
\addplot [color=black, mark=o, mark options={solid, black}, forget plot]
  table[row sep=crcr]{%
0.309016994374948	0.951056516295154\\
0.265115830725854	0.964216571263189\\
0.296851532460923	0.954923644945292\\
0.392718721303605	0.919658635547811\\
0.519296828849877	0.854593940738209\\
0.630793525981774	0.77595072496872\\
0.68144750810473	0.731866991807838\\
0.729428687256814	0.684056861822758\\
0.888386102650164	0.459097084087944\\
0.706514722799897	-0.707698344259039\\
0.936797351650631	0.349872436668516\\
0.963353495773664	-0.26823505024635\\
0.492888668840132	0.870092385973468\\
-0.922356956951175	-0.386338768393451\\
0.691235718559843	0.722629352702378\\
-0.629114989026446	-0.777312247801522\\
-0.266591891429652	0.963809505775888\\
0.604225736676516	-0.79681318961079\\
};
\addplot [color=black, dashed, mark=o, mark options={solid, black}, forget plot]
  table[row sep=crcr]{%
0.309016994374948	1.95105651629515\\
0.175320459892549	1.96017680710441\\
0.0310067596273049	1.91893949375538\\
-0.125349647149275	1.77499784639259\\
-0.270533358804914	1.46791952692323\\
-0.350391159922618	0.96902224577565\\
-0.225191669853146	0.309960128872585\\
0.332330239226554	-0.233719158005587\\
1.13637255651905	-0.509666415790914\\
1.20069978649088	0.161658385067158\\
0.508491940440797	-0.553761601384708\\
1.96315673827914	-0.248398813820785\\
0.503539201881182	-0.129850895491007\\
-0.130956663192619	-0.997637038502955\\
1.60992921285345	1.11760057132798\\
0.12248859706412	-0.117697171127405\\
0.0849270320498348	1.89999028235476\\
0.205820859750963	0.12039646692469\\
};
\addplot [color=black, draw=none, mark size=1.8pt, mark=square*, mark options={solid, fill=black, black}, forget plot]
  table[row sep=crcr]{%
0.309016994374948	0.951056516295154\\
};
\addplot [color=black, draw=none, mark size=1.8pt, mark=square*, mark options={solid, fill=black, black}, forget plot]
  table[row sep=crcr]{%
0.309016994374948	1.95105651629515\\
};
\node[right, align=left]
at (axis cs:0.3,1.5) {t=0.0};
\addplot [color=red, draw=none, mark size=2.5pt, mark=*, mark options={solid, fill=red, red}, forget plot]
  table[row sep=crcr]{%
0.888386102650164	0.459097084087944\\
};
\addplot [color=red, draw=none, mark size=2.5pt, mark=*, mark options={solid, fill=red, red}, forget plot]
  table[row sep=crcr]{%
1.13637255651905	-0.509666415790914\\
};
\node[right, align=left, font=\color{red}]
at (axis cs:1.1,0.4) {t=0.9};
\addplot [color=black, draw=none, mark size=3.5pt, mark=diamond*, mark options={solid, fill=black, black}, forget plot]
  table[row sep=crcr]{%
0.604225736676516	-0.79681318961079\\
};
\addplot [color=black, draw=none, mark size=3.5pt, mark=diamond*, mark options={solid, fill=black, black}, forget plot]
  table[row sep=crcr]{%
0.205820859750962	0.12039646692469\\
};
\node[right, align=left]
at (axis cs:0.5,-1.2) {t=2.0};
\end{axis}
\end{tikzpicture}%

%% file: figures/low_order_norms.tex
%
%
\begin{tikzpicture}

\begin{axis}[%
width=0.951\fwidth,
height=0.75\fwidth,
at={(0\fwidth,0\fwidth)},
scale only axis,
xmode=log,
xmin=0.007,
xmax=0.35,
xminorticks=true,
xlabel style={font=\color{white!15!black}},
xlabel={$h$},
ymode=log,
ymin=8,
ymax=1000,
yminorticks=true,
ylabel style={font=\color{white!15!black}},
ylabel={$\|\bm q^{n+1}\|_2$},
axis background/.style={fill=white},
legend style={legend cell align=left, align=left, draw=white!15!black}
]
\addplot [color=black, dotted]
  table[row sep=crcr]{%
0.00698936879834068	8.72566280013394\\
0.0181605506838155	8.73548683087083\\
0.0272582462253452	8.7579044473126\\
0.0360462528233301	8.79205581546224\\
0.0445627128586674	8.83704852869237\\
0.0528238886907962	8.89213604393359\\
0.0611084663068251	8.95910606951827\\
0.0691185460155218	9.03562361767026\\
0.0769858825855495	9.12288269825125\\
0.0846186911278534	9.22012555499615\\
0.0919606465130903	9.3267792735821\\
0.0991700372344537	9.445858034367\\
0.105971070265574	9.57346000066222\\
0.112353295507796	9.70932138755878\\
0.118309225786831	9.85323717015391\\
0.123834655040543	10.0050558234264\\
0.128928667224003	10.1646730596531\\
0.133593513445331	10.3320263362633\\
0.13783441517663	10.5070908049491\\
0.141659307017601	10.6898766404265\\
0.145078525588772	10.8804273585295\\
0.148104456890825	11.0788186825155\\
0.150833029051044	11.2921739665672\\
0.153173637122528	11.5141677562874\\
0.155144290765504	11.7449854289138\\
0.156810264219193	11.9927349523126\\
0.158158368824086	12.2586064327613\\
0.159178445576688	12.5439104473346\\
0.159863323797035	12.8500921211349\\
0.160208721033802	13.1787465601201\\
0.160208543901393	13.5414039631123\\
0.159836993574535	13.9414968872879\\
0.15907345223395	14.3828653153849\\
0.157872058049532	14.8812455021339\\
0.156199768242892	15.4436329074476\\
0.154032111714124	16.0782627567162\\
0.151300334819245	16.8086120964668\\
0.147978313570413	17.6494042219659\\
0.144051528147603	18.6187640394717\\
0.139521296158675	19.7380220933059\\
0.13434407717715	21.0482933762207\\
0.12853929518642	22.5864857998172\\
0.122114982109277	24.4071518093902\\
0.115162662792941	26.5585534254153\\
0.107682206967477	29.1395414106126\\
0.0998035266349816	32.2358618737839\\
0.0916748262699798	35.9517701075765\\
0.0833492604958582	40.4758609469967\\
0.0750164536190409	45.9795797189721\\
0.0667342650294566	52.7848116093424\\
0.0586054439979385	61.3097195418661\\
0.0507601663460157	72.1039615308383\\
0.0432662687697468	86.0495133539667\\
0.0361973092163493	104.477429396491\\
0.0296410389364306	129.407268848084\\
0.0236391448794215	164.33600066107\\
0.0182469232628454	215.295601472839\\
0.013507296848321	293.670527791713\\
0.00945191113347638	423.114578626383\\
0.00699911288883229	574.196017993515\\
};
\addlegendentry{\VarTay}

\addplot [color=black, dotted, forget plot]
  table[row sep=crcr]{%
0.350035579928447	33.8724412871051\\
0.308341122233988	37.4787577478995\\
0.270520631567201	41.7237694694116\\
0.235875649921422	46.8229680946781\\
0.203940827352903	53.0765794016175\\
0.174435412616944	60.9108463812068\\
0.147325105719732	70.8988711486497\\
0.122493812031178	83.9520296305476\\
0.0998781058852543	101.515825570557\\
0.0795239768177115	125.894434898177\\
0.0614217959445455	161.182737014637\\
0.0456051629997537	214.982148239435\\
0.0320948710462025	302.965439023559\\
0.0209550609137784	460.890374035838\\
0.0122055688439256	787.097786148978\\
0.0095899726454218	1000.19297523062\\
};
\addplot [color=blue, dashed]
  table[row sep=crcr]{%
0.00685731804267354	8.72731513756217\\
0.0150897654773787	8.73964235473037\\
0.0220921652688236	8.76382198055129\\
0.0286020436685504	8.79872351287103\\
0.0348396283552082	8.84472619837382\\
0.0406670782610733	8.90047367568751\\
0.0459782256333876	8.96400241934952\\
0.0509231372528873	9.03642664826395\\
0.0554668239016879	9.11707888388055\\
0.059609073368764	9.20570342217409\\
0.0633545951956752	9.30213837585326\\
0.0665852545188075	9.4019975308673\\
0.0694689715439238	9.50899659642805\\
0.0720181419663543	9.62320158717545\\
0.0741598453364949	9.73951642394609\\
0.0760197419673304	9.86269369402195\\
0.0775474127597503	9.98709344759093\\
0.0788414777591123	10.1181008213831\\
0.0799140685713308	10.2559114061541\\
0.0807773452739051	10.4007417901722\\
0.0814433984098386	10.5528288193979\\
0.0819241640767079	10.7124291926284\\
0.0822413762257305	10.8876177760333\\
0.0823842089651386	11.0798676974232\\
0.0823438743139714	11.2908360900767\\
0.0821011908966145	11.5315659418983\\
0.0816309084452043	11.8057898796959\\
0.0809122165341322	12.1178666618895\\
0.0798963443881988	12.484043274369\\
0.0785548580363733	12.9128967043002\\
0.0768232915897049	13.4279561879171\\
0.0746766447949111	14.0457218968402\\
0.0720101776004777	14.8145089745595\\
0.0687762204559903	15.7810881990427\\
0.0648362154918539	17.0471949890998\\
0.0600747325819203	18.7634703839854\\
0.0543076743060015	21.2226026195567\\
0.04716704469904	25.1113326536669\\
0.029316959903449	44.1714408448422\\
0.0277838774047992	47.305871818056\\
0.0269680985444317	49.3050971770254\\
0.0265779409488452	50.4892437726936\\
0.0264373823660085	51.1231865251014\\
0.0264313055039737	51.4173508670539\\
0.0264945030734899	51.5285075312912\\
0.0266299802869495	51.5140602445479\\
0.026882199701422	51.3222419800375\\
0.0273265742516496	50.8397470421171\\
0.0280783112838113	49.902337829443\\
0.0293114064797086	48.3062271488234\\
0.031329673635992	45.7920566974644\\
0.034833174891106	41.923662346149\\
0.0435935005938952	34.6543705058501\\
0.0552695914351018	28.3648679440998\\
0.0649293554213969	24.8319879346259\\
0.0744662241857399	22.2379879004365\\
0.0842313340194813	20.1950190085985\\
0.0943111521944882	18.5385834267686\\
0.104790877891274	17.1648668400815\\
0.115776926490227	16.0032013105747\\
0.127252619737681	15.015511512642\\
0.139342080747861	14.1623406938168\\
0.152174798553049	13.416066719393\\
0.165962153418338	12.7535234145954\\
0.180742037073333	12.1655449000016\\
0.196730282581684	11.6381817510084\\
0.214215111273573	11.1598605992088\\
0.233609664609106	10.7202371970495\\
0.255525794382982	10.3095739242057\\
0.280606591073093	9.92227555865846\\
0.309984490813327	9.54993524228142\\
0.344824816099674	9.1896204718024\\
0.350253012326757	9.13972757216346\\
};
\addlegendentry{\BE}

\addplot [color=red, dashdotted]
  table[row sep=crcr]{%
0.00697719752945109	8.72571513063209\\
0.0179426310173878	8.73584326310121\\
0.0268433457341043	8.75874522778724\\
0.0351630428062811	8.79245105461657\\
0.043195944928793	8.83668618482407\\
0.0512030944628969	8.89285111455546\\
0.0588794772047582	8.95896556942682\\
0.0660807456591693	9.03318260425946\\
0.0729615004358579	9.11677615015616\\
0.0794337850453251	9.20876668572027\\
0.0854509430103822	9.30835968442143\\
0.0910086528950632	9.4152591358354\\
0.0961083891905871	9.52924794728492\\
0.100757411138787	9.65019231646855\\
0.104967810203983	9.77803952032236\\
0.108755268401503	9.9128115535154\\
0.112137873001134	10.0545970861253\\
0.115135128603845	10.2035434781397\\
0.117767195094524	10.3598498163569\\
0.120054329347975	10.5237613883243\\
0.122016492168087	10.6955656920468\\
0.123673081553534	10.8755899343591\\
0.125086595987279	11.0710997512792\\
0.126212111553622	11.2762701520885\\
0.127093305923931	11.4992107662229\\
0.127721886887875	11.741647518992\\
0.128091248200512	12.0055454752106\\
0.128195377996838	12.3021289881802\\
0.128006597147808	12.6357377649941\\
0.127500161463719	13.0115909370036\\
0.126629153424721	13.4470677713631\\
0.12538400620152	13.9410003326346\\
0.123666423514663	14.5306354041682\\
0.121458670938567	15.2255642289285\\
0.118588710000927	16.0930170458905\\
0.114904394405019	17.2113617588685\\
0.105080924227	20.7470855478374\\
0.104185556612526	21.2679955362738\\
0.103881493517348	21.5889143312114\\
0.10393453003105	21.7814347526198\\
0.104209270542448	21.887704784047\\
0.104648872202044	21.9374837699153\\
0.1053502518359	21.9399411250045\\
0.106375873705717	21.8804187018329\\
0.107897701926068	21.734401909444\\
0.110177738929317	21.4610547137023\\
0.113568248175108	21.009501653253\\
0.11863003435102	20.315290654858\\
0.126605324119542	19.2615977027357\\
0.142357235833907	17.4387299579825\\
0.168580582121629	15.1177729248757\\
0.185564538715694	13.9821460543064\\
0.201195199161965	13.128813930202\\
0.216248510723084	12.4463396109031\\
0.230954983170681	11.8877438550464\\
0.245411255953649	11.4253173151569\\
0.259798215436706	11.0369459052293\\
0.27420781998837	10.7090599911344\\
0.288623765725655	10.4334898040658\\
0.303232831377739	10.2004008337001\\
0.317923229363902	10.0066036370066\\
0.33290672615305	9.84526931245903\\
0.348244822260573	9.71293555903776\\
0.350132098341198	9.69871104534957\\
};
\addlegendentry{\CN}

\addplot [color=black]
  table[row sep=crcr]{%
0.00473370336464369	8.72596605888736\\
0.0176850110116592	8.73502492028521\\
0.0268203397511629	8.75754298723497\\
0.0354788637070035	8.79149266073306\\
0.0439114820202223	8.83656187244497\\
0.0521649641056549	8.89248292785499\\
0.060247348321093	8.95911302986851\\
0.0681314912881166	9.03623590756032\\
0.0757857579283129	9.12368365935\\
0.0831535380412392	9.22099887904558\\
0.0902110102049329	9.3280817969699\\
0.0969194589854205	9.44457021195093\\
0.103256004356775	9.5702273977215\\
0.109212895427768	9.70501619397868\\
0.114788797842412	9.84898845068347\\
0.119982232357274	10.0021101022683\\
0.124811683974615	10.1648980120379\\
0.129291399784911	10.3378318440271\\
0.13344500579386	10.5218894756104\\
0.1372900680274	10.7180247467192\\
0.140847556031106	10.92760382832\\
0.14414575647203	11.1529004217601\\
0.147210287187408	11.396781130576\\
0.150066088509306	11.6629804437535\\
0.152746985004287	11.9573699000621\\
0.155296747320808	12.288982479082\\
0.157804900984103	12.6764184570451\\
0.160553587458022	13.177452655115\\
0.165693004422355	14.1624246378657\\
0.167500293097347	14.4214496786206\\
0.169173855972089	14.5983191540657\\
0.170841176937675	14.7202729409868\\
0.17257350867992	14.7998464358688\\
0.17443936729559	14.8433607995431\\
0.176514829576004	14.8525665361106\\
0.178893767057531	14.8254169108127\\
0.181684906424852	14.7565799375311\\
0.185020981803649	14.6379988640246\\
0.189068165868075	14.4592494668901\\
0.19403302207891	14.2081806546454\\
0.200164713544115	13.8722212866522\\
0.207684582869817	13.4439167194582\\
0.21611327007564	12.9578416364556\\
0.21884717161152	12.7948419033198\\
0.216360382681804	12.9424664323173\\
0.208962191934363	13.4435678470209\\
0.200533413101059	14.1091091597635\\
0.192712473027948	14.8346523077748\\
0.186380252344204	15.5275623681129\\
0.181904493748822	16.110819409897\\
0.179159688848127	16.5510562377411\\
0.177769850587005	16.8543785753623\\
0.177321543456564	17.0470236656479\\
0.177470358903059	17.1599825111378\\
0.178007375053089	17.2206717364347\\
0.17892636565716	17.2439194417379\\
0.18038613227071	17.2262397428219\\
0.182649023030969	17.1506251724941\\
0.186074545882517	16.9911230885008\\
0.19114354216497	16.7157098087649\\
0.1985667392156	16.2857940066265\\
0.209640003399173	15.6459426986492\\
0.229343931141094	14.5944316352286\\
0.26045648049396	13.2306725820908\\
0.279468568548391	12.5683122323569\\
0.297473020777285	12.0424147000032\\
0.315668147271269	11.5943896062222\\
0.334685129529596	11.1992904828762\\
0.350006838735574	10.9253742807652\\
};
\addlegendentry{\trbdf}

\end{axis}
\end{tikzpicture}%

%% file: figures/high_order_norms.tex
%
%
\begin{tikzpicture}

\begin{axis}[%
width=0.951\fwidth,
height=0.75\fwidth,
at={(0\fwidth,0\fwidth)},
scale only axis,
xmode=log,
xmin=0.007,
xmax=0.35,
xminorticks=true,
xlabel style={font=\color{white!15!black}},
xlabel={$h$},
ymode=log,
ymin=8,
ymax=1000,
yminorticks=true,
axis background/.style={fill=white},
legend style={legend cell align=left, align=left, draw=white!15!black}
]
\addplot [color=black, dotted]
  table[row sep=crcr]{%
0.00473370336464369	8.72605416131232\\
0.0179556298117626	8.73513342532619\\
0.0272459163419712	8.757750005587\\
0.0361156097823953	8.79187925596791\\
0.0448068231900564	8.83719742875475\\
0.0534046533811075	8.89365617091417\\
0.0619215149487312	8.96117294904987\\
0.0703612071177945	9.03981025180928\\
0.07869733447546	9.1294873633975\\
0.0869069339502662	9.23019440055067\\
0.0949852564646178	9.34220887154373\\
0.102911134223812	9.46564175576834\\
0.110678228046951	9.60081754399786\\
0.118297126500504	9.74839428766568\\
0.125791715575023	9.90939929782672\\
0.133185052424534	10.0849793409644\\
0.140518182647435	10.2768644500443\\
0.147838231495779	10.4871976194602\\
0.155209418614692	10.7189577252146\\
0.162696675095962	10.9756246093419\\
0.170361323232494	11.2611067445631\\
0.178227750209214	11.5785267495486\\
0.186251370962714	11.9287156870016\\
0.194257944580205	12.3070168370282\\
0.201944391594015	12.7020459296939\\
0.208936519489506	13.0967387416962\\
0.214941440074035	13.4750016981679\\
0.219802819591664	13.8248593154573\\
0.223538431605292	14.1424882708904\\
0.22628652762735	14.4320199038663\\
0.228221703681683	14.7025851074208\\
0.229500679221797	14.9667341047177\\
0.230225142327076	15.2386560380822\\
0.230434842774105	15.5341768427802\\
0.230109622674654	15.8728098744297\\
0.229169208975527	16.2789251157474\\
0.227479157310919	16.7830932261491\\
0.224865059255383	17.4229078792491\\
0.221137296068154	18.2434702734333\\
0.216103126646395	19.3028633784231\\
0.209475474738589	20.7038172369169\\
0.200126701466823	22.8044148510874\\
0.183049311807825	27.5192770337479\\
0.177275563954388	29.6180642567903\\
0.173245278475398	31.392862930757\\
0.170503548839842	32.88101384636\\
0.168749786626139	34.1121579357553\\
0.167757489686583	35.1204011443173\\
0.167357710760779	35.9381201695541\\
0.167423125980862	36.5945602416537\\
0.167859190449465	37.1167358188355\\
0.168601328587104	37.5279552163378\\
0.169616874437895	37.8480320089912\\
0.170908800282427	38.0916119103567\\
0.172513731859939	38.2661726008768\\
0.174503499805568	38.3718982856897\\
0.17697719860149	38.401362325677\\
0.180062678783292	38.3409105075431\\
0.183922565162712	38.171686010681\\
0.188758178971539	37.8708542851963\\
0.194829688032912	37.4123412588581\\
0.202481637835618	36.7673255993923\\
0.212175371244897	35.9053729321254\\
0.224556007234208	34.7946605653661\\
0.240547721197589	33.4032437458047\\
0.26151819357369	31.7006160192603\\
0.289475767658841	29.6660329989926\\
0.327089080816928	27.315172453554\\
0.350004043554996	26.0672314537383\\
};
\addlegendentry{\RadThree}

\addplot [color=blue, dashed]
  table[row sep=crcr]{%
0.00699690951150077	8.72566615326732\\
0.018200625243884	8.73556724893528\\
0.0274689491327388	8.75850392698725\\
0.0363211360830503	8.79293916397126\\
0.0449911761957111	8.83854394721055\\
0.0535633094636057	8.89526617863883\\
0.0620538704381994	8.96305581894311\\
0.0704549441369722	9.04188111045671\\
0.0787429057461951	9.13168423670993\\
0.0868899259604735	9.23242441138634\\
0.0948693011048509	9.3441110692039\\
0.10265572409241	9.46678800636044\\
0.110228709373312	9.60056611561891\\
0.117574062108719	9.74564769043148\\
0.124693911417426	9.90256941235277\\
0.131588659289915	10.0718823256492\\
0.138272028349572	10.2545187461754\\
0.144765318496935	10.4517613922257\\
0.151099702977871	10.6654186616985\\
0.157315945327931	10.8979623939659\\
0.163463394467893	11.1526314589057\\
0.169621271486301	11.4345359805531\\
0.175910136748283	11.751688790385\\
0.182556071826059	12.1189888212282\\
0.19017959482413	12.5762539274375\\
0.216550333149165	14.2273327177338\\
0.218134378054914	14.3803793766296\\
0.218962707736543	14.5228868803931\\
0.219255577502506	14.6769976814233\\
0.219028087254589	14.8606258393273\\
0.218175755989047	15.0951849827386\\
0.21648689771298	15.4105696287326\\
0.213633974873452	15.8510545243652\\
0.209173218204933	16.4821400785482\\
0.202612310632984	17.3942053416853\\
0.193630680632342	18.6918249039145\\
0.182351487102556	20.469285352664\\
0.169323029482964	22.8037271225292\\
0.155218737012366	25.7794563954513\\
0.140650937100798	29.5065665926611\\
0.126139207355064	34.1224185873143\\
0.112100266716703	39.7914953353096\\
0.098827808745603	46.7157847231733\\
0.0864895527768836	55.1575500637727\\
0.0763639847711613	64.2166957315757\\
};
\addlegendentry{\RadFive}

\addplot [color=blue, dashed, forget plot]
  table[row sep=crcr]{%
0.0763639847711613	64.2166957315757\\
0.0641505414110344	79.1036243558196\\
0.0548260097973315	95.089236580016\\
0.0463977375214931	115.228840776496\\
0.0388105120145753	141.029359824327\\
0.0320017698367056	174.829236464999\\
0.0259097166588422	220.432309565535\\
0.020461714222859	284.658259721252\\
0.0155562728119179	381.75265244761\\
0.0109838110263803	552.230039672792\\
0.00699998655489947	888.096753932263\\
};
\addplot [color=red, dashdotted]
  table[row sep=crcr]{%
0.00699592916817408	8.72566613085601\\
0.0181984536917894	8.73556341423686\\
0.0274484022664385	8.75843800299147\\
0.0362951509441933	8.79281939581446\\
0.0449577045169891	8.83834389176047\\
0.0535246910055264	8.8949841253749\\
0.0620107086916088	8.9626879014099\\
0.0704008926584577	9.04136332565491\\
0.0786835880687744	9.13107727863972\\
0.0868078824307014	9.23153497727585\\
0.094759766598516	9.34288599431847\\
0.102502420411827	9.46501135389538\\
0.110022635424791	9.59814484691982\\
0.117282559566725	9.74205389179914\\
0.124292160206186	9.8973982796118\\
0.131033651314278	10.0642982104142\\
0.137518391018507	10.243585604639\\
0.143744816413289	10.4358503798307\\
0.149712874772909	10.6417417959577\\
0.1554450482565	10.8628298993994\\
0.160950189229016	11.1004772019865\\
0.166245385114457	11.3566952630247\\
0.17133541658344	11.6332730694964\\
0.176231474855977	11.9326886821595\\
0.180949200831831	12.2582994716842\\
0.185493585192721	12.6134801241027\\
0.189873290483775	13.0026910803212\\
0.194091770965599	13.4310165003388\\
0.198142610720996	13.9036154923948\\
0.202033296405239	14.4287037730159\\
0.205748264175229	15.0133067674293\\
0.209286454252783	15.6685124393059\\
0.212642357360173	16.4077104792571\\
0.215824215812747	17.2515309025328\\
0.218871360727743	18.2366188400608\\
0.221986853365663	19.4744384668971\\
0.227929698734848	22.0909481573017\\
0.229622295146982	22.5843594688635\\
0.231166472189382	22.8825663168768\\
0.232651664763407	23.0559756352898\\
0.234134996619006	23.1430977853178\\
0.235698152084508	23.1642639450055\\
0.237448346708207	23.1235291912399\\
0.239522436587621	23.010565145431\\
0.242056428505332	22.804806640474\\
0.245197190914222	22.4792607711827\\
0.249097079305248	22.004667326084\\
0.253865983559311	21.3597552020523\\
0.259523137752296	20.5424629535252\\
0.265946027149232	19.5787090163122\\
0.272875017978858	18.5186822814549\\
0.279996035242299	17.4189321724129\\
0.287029646627801	16.3252843967929\\
0.293763188578462	15.2678391556505\\
0.300036680753612	14.2647755092806\\
0.305726987675885	13.3264107153891\\
0.310717745351485	12.4611183829196\\
0.314881946476183	11.679414646101\\
0.318086783727075	10.9951611889846\\
0.320207815091812	10.4279396046687\\
0.321214371513393	9.99413917873586\\
0.321244725498234	9.69456124253001\\
0.320578749627601	9.50926699431556\\
0.319517160740916	9.40798713369503\\
0.318281470544999	9.36202226890122\\
0.316928082961288	9.35082700541727\\
0.315320572330103	9.36775098009363\\
0.313271443781656	9.41993593277378\\
0.310627477043617	9.52279504065879\\
0.307346356018497	9.69382523453607\\
0.303536433230132	9.94646092958726\\
0.299350227888253	10.2922768763311\\
0.29480884102283	10.755148376356\\
0.289679232980146	11.3956899088512\\
0.283187561552601	12.3842331229345\\
0.256471889222387	17.9914706011841\\
0.248432162619269	20.0085979328295\\
0.240584776898379	22.0969279414561\\
0.232741044808749	24.3040450934278\\
0.224814539154933	26.6637777061245\\
0.216763471784488	29.2061709337063\\
0.208562933628534	31.9637853296423\\
0.200203596945872	34.9716896589022\\
0.191683500632818	38.2704723701979\\
0.183007153944554	41.9075497710744\\
0.174184311439968	45.9390995113689\\
0.165231722976833	50.4311181836346\\
0.156167401065562	55.4647850974133\\
0.147015876909776	61.1374009931063\\
0.137803213657532	67.5698598314798\\
0.128561619861825	74.9099708787892\\
0.119324805639454	83.3440138617263\\
0.110131160033826	93.1057215590578\\
0.101020742176853	104.494796656267\\
0.0920374269765998	117.89671198435\\
0.0832262217720298	133.817557460906\\
0.0746338175182332	152.929904076355\\
0.0663080172413331	176.142452275855\\
0.0582965423800148	204.708123065287\\
0.0506464086680803	240.392331124813\\
0.0434033296330198	285.744852119205\\
0.0366107985876476	344.558049239234\\
0.0303093703796034	422.668746430737\\
0.0274329235720058	470.24668856995\\
};
\addlegendentry{\GauLegFour}

\addplot [color=red, dashdotted, forget plot]
  table[row sep=crcr]{%
0.0274329235720058	470.24668856995\\
0.0218077580434419	599.549474068623\\
0.0168927797280135	783.004493846566\\
0.0133371142285949	1000.00472767134\\
};
\addplot [color=black]
  table[row sep=crcr]{%
0.00699753810773066	8.72566616759373\\
0.0182099290723165	8.73558355417561\\
0.0274772854167825	8.75853062594431\\
0.0363183696308939	8.79292669186804\\
0.044994308720021	8.83856339447566\\
0.0535623541413135	8.89526182446346\\
0.0620507321560374	8.96303536540533\\
0.070450924738656	9.04185596557162\\
0.0787360626168857	9.131637113638\\
0.0868891117879584	9.23247654912062\\
0.0948676865734081	9.34420431833606\\
0.102646724003073	9.46684920493991\\
0.110207212531692	9.60053609771457\\
0.11752486142563	9.74526082735724\\
0.124608804755746	9.90167245509717\\
0.131459695782454	10.0703475240546\\
0.138070944607387	10.2516933062829\\
0.14447012808286	10.4471283755396\\
0.150670915455371	10.6579095974483\\
0.156689315815591	10.8855454635253\\
0.162564241048286	11.1327462294916\\
0.168322073448537	11.4022206385951\\
0.174001501338419	11.6977516806063\\
0.179663339174505	12.025093470526\\
0.185388345700679	12.3923913728765\\
0.19132009676941	12.813650977098\\
0.197768030799146	13.3179636459123\\
0.206057318807075	14.0233452043218\\
0.216886889087988	14.9509310292417\\
0.221495525658258	15.2917885129513\\
0.225136410608761	15.5110301033809\\
0.228236243488226	15.6503134166655\\
0.230989542974031	15.7298687145926\\
0.233517452991049	15.7615389053975\\
0.235917264396065	15.7518695875364\\
0.238273384900377	15.702721393298\\
0.240654661194748	15.6118925014136\\
0.243120029701852	15.4736709165911\\
0.245707716754371	15.280262919847\\
0.248443362894808	15.0226025490137\\
0.251336335730131	14.6917712100332\\
0.254372388949674	14.2812951541796\\
0.257508462148285	13.7895907083558\\
0.260660485944909	13.223306962065\\
0.263655581589568	12.6072169641885\\
0.266139132863403	12.0064389119074\\
0.267558079172372	11.5478489201803\\
0.267773484273769	11.311494152032\\
0.267323563118011	11.2324168942965\\
0.266650806503394	11.2254250812379\\
0.265611667424723	11.2655643381399\\
0.263793323635227	11.3876584578734\\
0.260518413952815	11.6735962699344\\
0.254634837880091	12.2864144662035\\
0.245476725541122	13.4159521161947\\
0.236485435832685	14.7558207831979\\
0.23119418952204	15.7198725912855\\
0.228807317348028	16.2831469655602\\
0.227693219081797	16.6866375034187\\
0.227159666560153	17.0954879526785\\
0.226964801531246	17.7351400734504\\
0.226518706839558	21.1901933131008\\
0.225611896673145	23.3608939379334\\
0.224227645636137	25.6877785296825\\
0.222404407832017	28.1419280376853\\
0.220174423700696	30.7122893181563\\
0.217560833734621	33.4008517763746\\
0.214578490011466	36.2186632522182\\
0.211237361083122	39.1819274970817\\
0.207542807220337	42.3122769130208\\
0.203497935072009	45.635530814349\\
0.199104046684264	49.1824541466726\\
0.194361444668124	52.9893870859035\\
0.189270487030202	57.0989270127945\\
0.183831817275933	61.5615605015593\\
0.178047473026394	66.4369436626141\\
0.171921278452194	71.7962134901948\\
0.165460612674827	77.7236327709446\\
0.158675397791955	84.3213854213591\\
0.15158066627681	91.7120214703316\\
0.144197028130169	100.043685674062\\
0.136550729887057	109.497230491708\\
0.128674619413273	120.294115202092\\
0.120609027100511	132.706273788065\\
0.112400715694603	147.071633071381\\
0.104102712278986	163.812992644588\\
0.0957738915217132	183.46305029223\\
0.0874770570076614	206.701505519322\\
0.0792772343984785	234.406415786434\\
0.071240232352221	267.726042128234\\
0.0634304054320435	308.186559892492\\
0.0559087632202083	357.8526508501\\
0.0522687156249924	387.031558883822\\
};
\addlegendentry{\GauLegSix}

\addplot [color=black, forget plot]
  table[row sep=crcr]{%
0.0522687156249924	387.031558883822\\
0.0449559070104996	459.971161996295\\
0.0384087145729025	548.906090072165\\
0.0323187036564123	664.064939673237\\
0.0267182015623641	816.414463207152\\
0.0221049591321448	1000.00425271711\\
};
\end{axis}
\end{tikzpicture}%

%% file: figures/low_order_traj_errors.tex
%
%
\begin{tikzpicture}

\begin{axis}[%
width=0.951\fwidth,
height=0.75\fwidth,
at={(0\fwidth,0\fwidth)},
scale only axis,
xmode=log,
xmin=0.007,
xmax=0.35,
xminorticks=true,
xlabel style={font=\color{white!15!black}},
xlabel={$h$},
ymode=log,
ymin=1e-06,
ymax=1000,
yminorticks=true,
ylabel style={font=\color{white!15!black}},
ylabel={Relative error},
axis background/.style={fill=white},
legend style={at={(0.03,0.5)}, anchor=west, legend cell align=left, align=left, draw=white!15!black}
]
\addplot [color=black, dotted]
  table[row sep=crcr]{%
0.00698936879834067	1.72438587670581e-06\\
0.00896229975504268	3.52174355969678e-06\\
0.0112061681347665	6.64731825122414e-06\\
0.013718641596227	1.17360638307378e-05\\
0.0164980430295553	1.95909261603393e-05\\
0.0198198309398563	3.24081965816543e-05\\
0.0239568616140422	5.41803641970116e-05\\
0.0300059331450632	9.90948297174839e-05\\
0.0409894662366299	0.000228655055639082\\
0.0462131292402196	0.000317410678929674\\
0.0508941813320803	0.000415786380056485\\
0.0553068633412254	0.000528108362538943\\
0.0594497852381056	0.000654267854185946\\
0.0635979495369462	0.00080463726413662\\
0.067745247097764	0.000983847614482325\\
0.0718546000082679	0.00119575213423752\\
0.075913687986323	0.00144565651153624\\
0.0799109356309032	0.00173954588737309\\
0.084100649904167	0.00210964451380261\\
0.0882066631275793	0.00254783582135722\\
0.0924537335182181	0.00309818185825787\\
0.0965733788482757	0.00374888596554893\\
0.100793896812318	0.00456434310246781\\
0.1050882632414	0.00558802832032441\\
0.1094275932872	0.0068742473731255\\
0.113781350939696	0.00849024108208342\\
0.118309225786831	0.0106195494603064\\
0.122764099879368	0.0133018230424195\\
0.127110879216908	0.0166681298237562\\
0.131469726932872	0.0210503224677514\\
0.135624935655011	0.0265203811223737\\
0.139670558071614	0.0335547408860125\\
0.143418895028297	0.0422408451962316\\
0.146840639418851	0.052868219448576\\
0.149910137129037	0.0657473675267701\\
0.152605763806681	0.0812051233364467\\
0.154910151343037	0.0995796697094428\\
0.156810264219193	0.121215838081804\\
0.158263080534629	0.145737845855576\\
0.159321443318197	0.174041722959735\\
0.159970576162066	0.20556304669263\\
0.160243754994118	0.241345220186883\\
0.160150057816131	0.280630742737514\\
0.159704335884862	0.324638174006916\\
0.158935281010949	0.372412024107758\\
0.157841173917179	0.425350132491373\\
0.156466983605796	0.482303717415236\\
0.154796275716661	0.544879585781357\\
0.152887027305464	0.611745488327848\\
0.150709858975418	0.68474596630064\\
0.148332874953501	0.762380266993491\\
0.145715407891894	0.846770094693599\\
0.142933135455485	0.936267523404471\\
0.139939434446372	1.03320035880796\\
0.136751498605731	1.13791038239393\\
0.133381331942569	1.25099102311861\\
0.12991697354688	1.37046875300386\\
0.126296798030757	1.49952555868086\\
0.122498452814868	1.6403253681673\\
0.118557719052386	1.7931202284172\\
0.114508762555259	1.95825732311911\\
0.110377893149128	2.13643395639855\\
0.106123561596883	2.33160836312212\\
0.101778995733112	2.54488937353492\\
0.0973153476838125	2.78088900251828\\
0.0927731697809178	3.04130809711831\\
0.0881915283418081	3.32806280107653\\
0.0835555544332901	3.64703529319961\\
0.0788607740409789	4.00494858006411\\
0.0741126933594556	4.40955359052474\\
0.0693260243957882	4.86986095244872\\
0.0645235412344015	5.39642627592357\\
0.0597345900352251	6.00170002229509\\
0.0549933131473637	6.70045536819732\\
0.050280545432021	7.52096104346058\\
0.045652685276076	8.48659543711423\\
0.0410973193150446	9.64418295498291\\
0.0366497644278462	11.0462823730653\\
0.03232326651238	12.7742715896832\\
0.0281555693886232	14.9343499267646\\
0.0241531830065471	17.70320835537\\
0.0203605321998668	21.3236728516018\\
0.0167973149569893	26.2062419193245\\
0.0134963187843619	33.0209066623014\\
0.0104884812879425	42.9555251488807\\
0.0078043675652936	58.2760542306881\\
0.00699911288883228	65.1615988383127\\
};
\addlegendentry{\VarTay}

\addplot [color=black, dotted, forget plot]
  table[row sep=crcr]{%
0.350035579928447	3.64809823891989\\
0.340110646061333	3.71250376299902\\
0.330552810689066	3.76747498450142\\
0.32134017072815	3.81315597710754\\
0.312453324177734	3.84946187775651\\
0.304085762801319	3.87561851332024\\
0.296200965919711	3.89211131385108\\
0.288568581889698	3.8996508624077\\
0.280986080806841	3.8983570411931\\
0.27309171548541	3.88786257188348\\
0.262982188880393	3.86408560487027\\
0.252313875422863	3.83966771267237\\
0.247176025532178	3.8366523951573\\
0.243124224169359	3.84269765151881\\
0.239466969858495	3.85705537761579\\
0.236030442958015	3.88010270832702\\
0.232652977998177	3.91321444640576\\
0.229183812389273	3.95911101747264\\
0.225483095233689	4.02203971960733\\
0.221422971956115	4.10774488864832\\
0.216749685268721	4.22697349561899\\
0.211106907704197	4.39718722334719\\
0.202516967673121	4.7010458275712\\
0.186558679939266	5.3654251853455\\
0.177705775699793	5.77545885157938\\
0.169247604469568	6.19296960826042\\
0.160778207953758	6.63888279676793\\
0.152079359370992	7.13102150653114\\
0.14300318474064	7.68940222631172\\
0.133473187449354	8.33706499734039\\
0.123413266087646	9.10720794885589\\
0.112717387549083	10.0524254515544\\
0.101318928639087	11.2516901020053\\
0.0891600140103406	12.8363286797926\\
0.076163563735502	15.0488711841624\\
0.0620612864318583	18.4417591743972\\
0.045836290858435	24.8300360649912\\
0.0145644971119101	76.3171303547224\\
0.0083978968008478	131.45862440793\\
0.00699865345886625	157.481585235481\\
};
\addplot [color=blue, dashed]
  table[row sep=crcr]{%
0.00685731804267354	0.000413324709943259\\
0.0096922270847041	0.000845067177589161\\
0.0136850261981599	0.00173299287512397\\
0.0190614180930137	0.00347267714284461\\
0.0246528022616893	0.00598535624060144\\
0.0297612367832148	0.00895977504598369\\
0.0343801285819588	0.0122687911861167\\
0.0384604106720697	0.0157483815161063\\
0.0421849637815981	0.0194522405285927\\
0.0455654116924551	0.0233247673457928\\
0.0487966201402063	0.0275646228902398\\
0.0518640511311605	0.0321717629337493\\
0.0547668654755395	0.0371623328250845\\
0.0575050230000351	0.0425540061913525\\
0.0600793426251543	0.0483654361632744\\
0.0626368543702024	0.0550213887817186\\
0.0650140022072535	0.0621957047082392\\
0.0672142514526063	0.0699108106295946\\
0.0693555870209773	0.0786931271402983\\
0.0713084060403767	0.0881300858610108\\
0.0731715168339611	0.0988264622765835\\
0.0748380559958703	0.110301909642531\\
0.0763880782679052	0.123247608972525\\
0.0777997186529711	0.137822394756927\\
0.0790005100996648	0.153427463481184\\
0.0800438426342724	0.170905359631361\\
0.0808795158890626	0.18956276664344\\
0.0815442455633596	0.210353474579577\\
0.0820088060153008	0.232481707559535\\
0.0822945509087835	0.257023422521681\\
0.0823917077032178	0.284192311521913\\
0.0823000909507484	0.313032442819152\\
0.0820210065287285	0.344838521932546\\
0.081570723139	0.378532357663349\\
0.0809386317420933	0.41556787445986\\
0.0801552236466254	0.454742625341947\\
0.0791995111341547	0.497694591209149\\
0.0781138123599254	0.54309273737548\\
0.0768682533758954	0.592785664072166\\
0.075466225858071	0.647172047931003\\
0.0739125633215664	0.706701437474064\\
0.0722779573575974	0.76939711964882\\
0.0705013274216348	0.838388756418837\\
0.0686148437031144	0.913316762696622\\
0.0666341954421191	0.994515896133666\\
0.0645216698298538	1.084701935856\\
0.0622955209312868	1.18455524772755\\
0.0599752714781093	1.294848057792\\
0.0575809785924284	1.41645964158597\\
0.0551325295302822	1.55039769992083\\
0.0526048316216154	1.70056733803204\\
0.0499797886193186	1.87159555715779\\
0.0472859932136538	2.06610382064454\\
0.0445484728100299	2.28751308228114\\
0.0417497203882828	2.54411081164224\\
0.0388634619671472	2.84862821668618\\
0.0358399217592635	3.22335716918995\\
0.0324084535962997	3.74235387951392\\
0.0278138671094391	4.69345053628263\\
0.0269680985444317	4.93392714643745\\
0.0265779409488452	5.07137228606854\\
0.0264373823660085	5.14472925329366\\
0.0264313055039737	5.17857992440818\\
0.0264945030734899	5.19116425162151\\
0.0266204363170976	5.18961665323745\\
0.0268490199710201	5.1698561051343\\
0.0272447836093461	5.12057059993604\\
0.0278990240273117	5.02608219025479\\
0.0289221201152284	4.87067966637692\\
0.0304547792275782	4.64136056085172\\
0.0326090406166184	4.340401091554\\
0.0354640692521719	3.98443047483195\\
0.0389521686716768	3.60917590462019\\
0.0429467919833402	3.24543649133541\\
0.0472777520843065	2.91362137057455\\
0.0518729876048867	2.61664665974295\\
0.0566306870900771	2.35547685653388\\
0.0614321741683591	2.12924110175325\\
0.0662841007232369	1.93079716634647\\
0.0711652739447368	1.75572071589548\\
0.0760655551077058	1.60010071869369\\
0.0809225636870308	1.46223690352976\\
0.0857381367144909	1.33894564305875\\
0.0904547150229914	1.22906916300146\\
0.095089972025567	1.13000498810442\\
0.0996772474788713	1.03943257842022\\
0.104178619430263	0.956798796533216\\
0.108646170818564	0.880106614597549\\
0.112956414723256	0.810517502953679\\
0.117168064623904	0.746204398025825\\
0.12135616360137	0.685456608314981\\
0.125395016448465	0.629578526108884\\
0.129366934451179	0.576955653911295\\
0.133253700668881	0.527488152328887\\
0.137036822114572	0.481102677700092\\
0.140696006529038	0.43777119671213\\
0.14420549067674	0.397549854114107\\
0.14767723730384	0.359004023028817\\
0.151103732996027	0.322164843664787\\
0.15447741274196	0.287100511610054\\
0.157790687148639	0.253929058958913\\
0.161178648080874	0.221507208756756\\
0.164642984060159	0.190304061438734\\
0.168783559472124	0.15659679138024\\
0.174267984359988	0.121795872482251\\
0.176137316978602	0.113889890006314\\
0.177553081909513	0.109676561701481\\
0.178662464192706	0.107555706299988\\
0.179459320307164	0.106703064455219\\
0.180099483317521	0.1064309129478\\
0.180742037073333	0.10652739163208\\
0.181386991053476	0.10699228295918\\
0.182196573390063	0.1080849060516\\
0.183173068811124	0.110124180282408\\
0.18448359571149	0.114007465094614\\
0.186135578959825	0.120545585924606\\
0.18864267855149	0.133258483801332\\
0.199923494595918	0.209928699992837\\
0.203900705172342	0.238697748342515\\
0.207588996579915	0.26453647273628\\
0.211158501504977	0.288165857007456\\
0.214600497790874	0.309186498557831\\
0.217906687230111	0.327352458300106\\
0.221069744717831	0.342564118438067\\
0.224081541424704	0.354852563768376\\
0.226728953322006	0.363759119857947\\
0.229204035611275	0.370422614571087\\
0.231499806698621	0.375151732720981\\
0.233609664609106	0.378278635315347\\
0.23574163127418	0.380284757635646\\
0.237896113027264	0.381185215737561\\
0.240073535461518	0.381013565508092\\
0.242274344818936	0.379820307609866\\
0.244722806738734	0.377405979921279\\
0.247200790894413	0.373937932296378\\
0.250168331495707	0.368639545266958\\
0.25341274635007	0.361702518142949\\
0.257422473168957	0.3519273472669\\
0.262242330637791	0.339100904982966\\
0.26868306554553	0.32131336139757\\
0.292442188450962	0.264932081213991\\
0.298205285988147	0.255344227644866\\
0.303247073002005	0.248508560878643\\
0.307813442795611	0.243592369731556\\
0.311866607042601	0.240249361193791\\
0.315370079825293	0.238131239632967\\
0.318617297557551	0.236802349953955\\
0.321927136790093	0.236069460101485\\
0.325303163978305	0.235960011258343\\
0.328749315306981	0.236502268560062\\
0.331914401791623	0.2375717405762\\
0.335506205556011	0.239434280999119\\
0.33917868205539	0.2420351183912\\
0.3433026072356	0.245771628064612\\
0.347907047465948	0.250930183841275\\
0.350253012326757	0.253946404263693\\
};
\addlegendentry{\BE}

\addplot [color=red, dashdotted]
  table[row sep=crcr]{%
0.00697719752945108	1.06751372551411e-05\\
0.0086435754021289	1.96003813014783e-05\\
0.0103004482469899	3.20427013153683e-05\\
0.0122232809225397	5.1445941065244e-05\\
0.0144096300443247	8.05675663690165e-05\\
0.016857778350048	0.000122738846374032\\
0.0195668190446591	0.000181858677690578\\
0.0225366384124974	0.000262437939805547\\
0.0263056483919969	0.000389665820546573\\
0.0407880146049407	0.00118902922227198\\
0.044265462972516	0.00147780452679999\\
0.0474709355626961	0.0017903504170846\\
0.0506705300877456	0.0021553066796302\\
0.0538621515470569	0.00258202380672504\\
0.0570396558474933	0.00308047986910886\\
0.0604451166899381	0.00371123503762765\\
0.063797618539364	0.0044475078941855\\
0.0673347087859083	0.00537133004463461\\
0.0710316254433418	0.006528671451151\\
0.0748567243580847	0.00797326925601287\\
0.0787619716052566	0.00976125650269254\\
0.0829282912798572	0.0120954025500314\\
0.0870857920572408	0.0149688685063446\\
0.0911986453012727	0.0184842585886048\\
0.0952310575859309	0.0227574922089765\\
0.0991481155065139	0.0279183693737165\\
0.102916527886929	0.0341109029374443\\
0.106505220656512	0.0414934303980518\\
0.109885769506462	0.0502385128675603\\
0.113032678216943	0.0605326615024992\\
0.115923523810979	0.0725759870527305\\
0.118538994592672	0.0865818650890218\\
0.120862846369436	0.102776755909516\\
0.122881798409737	0.121400329427182\\
0.124537730597358	0.142007578882116\\
0.125892483573465	0.165425034806921\\
0.126939164204697	0.191925448868286\\
0.127656176823047	0.220884539677966\\
0.128076090036817	0.253372406053479\\
0.12819634356682	0.289720374696584\\
0.128024739317801	0.329130798060635\\
0.127567554072173	0.373004653217106\\
0.126847758418712	0.420434270929322\\
0.12588516210441	0.471635621160911\\
0.124660881145854	0.528482859196394\\
0.123213611033573	0.590014068469066\\
0.121506552183707	0.65868806749283\\
0.119580042655912	0.733863837361171\\
0.117423124700302	0.81721034477034\\
0.115032667593473	0.910414183336785\\
0.112338658544673	1.01857769958195\\
0.109087816606774	1.15734575246769\\
0.104667945414889	1.38445092609898\\
0.104006637684241	1.43819066848963\\
0.103869123246819	1.46800093467035\\
0.104006951072369	1.48332942523181\\
0.104279397835748	1.49030582063203\\
0.104648872202044	1.49244939190685\\
0.105187120579443	1.49034606776238\\
0.105961782836842	1.48252378474704\\
0.107111335131109	1.46599601782831\\
0.108744165632116	1.43769597384064\\
0.110960366169223	1.39512253522919\\
0.113765560774596	1.33840811583716\\
0.117089045697933	1.27007715595795\\
0.120770507983566	1.19480584544408\\
0.124654226631507	1.11690159134201\\
0.128660412346562	1.0386517664375\\
0.13275296317524	0.961118979044388\\
0.136772904190742	0.887347039399155\\
0.140801835077576	0.815717362127437\\
0.14475838806568	0.747520724091941\\
0.148636883120881	0.682621350402708\\
0.152430236694809	0.620907181545502\\
0.15604146040098	0.563680411791159\\
0.159547415784272	0.509467284075563\\
0.162939399778516	0.458217962324778\\
0.166208596093455	0.409898642595788\\
0.169346198959117	0.364489989167364\\
0.172441016086551	0.320617002612166\\
0.175390096645539	0.279665053701901\\
0.178184967771602	0.241652688214207\\
0.180817463629607	0.206609453077645\\
0.183279782703981	0.174576956763185\\
0.185669024725013	0.144308127257541\\
0.187876058859957	0.117253507264722\\
0.19000094016994	0.0923704778468048\\
0.192148455800719	0.0691384361776338\\
0.19607181241979	0.0402015172136342\\
0.196622668418391	0.039130823411932\\
0.196953883588315	0.0389958083444352\\
0.197174986549615	0.0391198554105872\\
0.197617896447887	0.0398699260517483\\
0.198284025044865	0.042157184506043\\
0.199398959701744	0.0484650925283182\\
0.20380414762952	0.0866079662324616\\
0.206098823361701	0.108569102739889\\
0.208534596487751	0.131603122931551\\
0.211115667155485	0.155139064588746\\
0.213727099282554	0.177681425270059\\
0.216490119166404	0.199856356817508\\
0.219287212675636	0.220311961731815\\
0.222118863768283	0.238779693330995\\
0.224860204619619	0.254369196783569\\
0.227507298324875	0.267203852578368\\
0.230056309090135	0.277476721237276\\
0.232503508601739	0.28543359400549\\
0.234845433349165	0.291354318370311\\
0.23694710383553	0.29532481278124\\
0.23906789559485	0.29812462825993\\
0.241073771661312	0.299740947828763\\
0.24296155546445	0.300425411915402\\
0.244864663324092	0.300374021022808\\
0.246920982489083	0.299575147711844\\
0.249134428352494	0.297966713857021\\
0.251790303737699	0.295183292613799\\
0.254903606769408	0.290992160174312\\
0.258637252837178	0.285036660336635\\
0.263612032895951	0.276222172880306\\
0.272355722135105	0.260377444639542\\
0.282257654437091	0.244449134456705\\
0.288955519346803	0.235607996601692\\
0.295021426133476	0.229005344101527\\
0.300926128949106	0.223750254627682\\
0.306838492247096	0.219476387561855\\
0.313139107909489	0.215785692873017\\
0.320663023401894	0.212203593442106\\
0.335896739167585	0.206143685690519\\
0.345001845590188	0.202373587106602\\
0.350132098341198	0.199989239895762\\
};
\addlegendentry{\CN}

\addplot [color=black]
  table[row sep=crcr]{%
0.00473370336464369	3.44753042166296e-06\\
0.00750020490266158	6.43401651636685e-06\\
0.00787532461798416	7.39343100213039e-06\\
0.0110806591766687	1.93451945729866e-05\\
0.0155977023243896	4.94347884896608e-05\\
0.0192314867125685	8.64782738595692e-05\\
0.0218510058234299	0.000120841032124716\\
0.0247981065578027	0.00016753055197096\\
0.0303177245136323	0.000279100576881189\\
0.0425368918573152	0.000657375044209036\\
0.0465947538705711	0.000834523339779394\\
0.050243434099841	0.00102308235216132\\
0.0537223440955048	0.00123381195763991\\
0.0571769782071933	0.00147855114812168\\
0.0606713715768381	0.00176854145945806\\
0.0643209745108406	0.00212503390965358\\
0.0680491893386463	0.00255567366935607\\
0.0720152476063089	0.00310058319084333\\
0.0762026586222338	0.00379058630890326\\
0.08059880942959	0.00466562742697788\\
0.0852406917061359	0.00579027732785086\\
0.0900960049265919	0.00723366487655868\\
0.0951400165190455	0.00908742829458639\\
0.100284708429843	0.011439526380509\\
0.105500335195409	0.014421437860857\\
0.110705789552229	0.0181598004274512\\
0.115872328689664	0.0228387146916101\\
0.120953213266928	0.028664418523743\\
0.12589842952878	0.0358694936513015\\
0.130683940898893	0.0447628192170947\\
0.135283317327746	0.0557135966383001\\
0.139677010780536	0.0691863729667768\\
0.143847267633934	0.0857574027397793\\
0.147794596013605	0.106234507211287\\
0.151549851200395	0.131905389797531\\
0.155236205270353	0.165535592171838\\
0.159779865049483	0.222606437644554\\
0.163010535062802	0.270987817842049\\
0.164733693257239	0.29423610060721\\
0.166105935921023	0.308821004130713\\
0.167281688491271	0.317806415230062\\
0.16830539476911	0.322914199750504\\
0.169190075142242	0.325398659864152\\
0.169958257863025	0.326222433729438\\
0.170691916446961	0.325956357294808\\
0.171493674417844	0.324603471741933\\
0.17243744164538	0.321746494907566\\
0.173573963381057	0.31672967864927\\
0.17493600380214	0.308804938039815\\
0.176537481385965	0.297311200229149\\
0.178373638047692	0.281846584088716\\
0.180411514898319	0.262488012491564\\
0.182591981622141	0.239862336410906\\
0.18484875364653	0.214926700175622\\
0.187116678139323	0.188769851432964\\
0.189338400114093	0.162434511823908\\
0.191469499090005	0.136787548658526\\
0.193478517435138	0.112488267653723\\
0.195348913968217	0.0899676398888503\\
0.197072103174323	0.0695473077458283\\
0.198667034586242	0.0513053661164796\\
0.200250092100735	0.0347624859099718\\
0.202150910502267	0.0220191328540731\\
0.202510929273439	0.0215938150435383\\
0.202700279080697	0.021712102930408\\
0.203023339266742	0.0224424747274817\\
0.203615004294027	0.0252756375484797\\
0.207731887149085	0.0645587629927732\\
0.209618389470334	0.0847006722233356\\
0.211694665472027	0.106626752565279\\
0.213983339679395	0.130198312171451\\
0.216543415502728	0.155801780583575\\
0.218809740035075	0.180805043365791\\
0.21849934400721	0.179484769325249\\
0.217339374878991	0.171512796214239\\
0.215324356817709	0.155250031906962\\
0.213055879415679	0.135198727659583\\
0.21070344934032	0.113424196207183\\
0.208082618259494	0.0895688536486089\\
0.204699688037775	0.065847977468759\\
0.203744335302214	0.0627234797063753\\
0.203153528026152	0.0620046989162375\\
0.20279766060172	0.0620618258568491\\
0.202391326696089	0.0625905647782619\\
0.201810467901863	0.0642019449483059\\
0.200973375307098	0.0682115435734616\\
0.199654574231506	0.0780181025408328\\
0.193266100531513	0.158533098209222\\
0.190585907972917	0.200619966229221\\
0.187771536816385	0.248321662805441\\
0.184720514045944	0.304115632507267\\
0.18113320360015	0.376642609016701\\
0.177541787632994	0.469495337024987\\
0.177325866286553	0.487142224010836\\
0.177593783229808	0.492306099982989\\
0.17794426349091	0.49247302354473\\
0.178475562725126	0.489933864632037\\
0.179416094174612	0.482342029465297\\
0.18096705342695	0.466257954757856\\
0.183180466648048	0.440022315202108\\
0.185904814224004	0.40543609473386\\
0.188926596516256	0.365929690184956\\
0.192101148879499	0.324281834372502\\
0.195365092521598	0.282184547082549\\
0.198727188326478	0.240453136304596\\
0.202325247461433	0.198785947763802\\
0.207334145096649	0.149228401387907\\
0.210941086242333	0.123252620298385\\
0.212867704884838	0.114042848152411\\
0.21428788310276	0.109649788489955\\
0.215388160467699	0.107678779309076\\
0.216230318538328	0.106993857751382\\
0.216891494327061	0.106932812510355\\
0.217568194305659	0.107278654890646\\
0.21840631294316	0.108235447611402\\
0.219511374221607	0.110288758159731\\
0.221031399364558	0.114316487929476\\
0.223426994838504	0.122604119968076\\
0.230174040429124	0.149443818849493\\
0.233201779696945	0.160203922717264\\
0.236007091803221	0.168726548052605\\
0.238700131632835	0.175513547222465\\
0.241341212749141	0.1809042840019\\
0.243978520723509	0.185177091639761\\
0.246687159268865	0.188602513942813\\
0.24958327252891	0.191427947044375\\
0.25292504637366	0.193942702379834\\
0.257609940187767	0.196753561931164\\
0.267039860816955	0.2022708061923\\
0.272667761200102	0.20625924679572\\
0.279535669830479	0.211845065241253\\
0.296679470843927	0.226149098557437\\
0.302607330223706	0.230154415626893\\
0.307794683176642	0.232947615309991\\
0.312528378393978	0.234823273253862\\
0.316950765608079	0.235941302715364\\
0.321166018548112	0.236403287090335\\
0.325249859463555	0.236267211420078\\
0.329281025406643	0.235554877081807\\
0.333328010168789	0.23425213035854\\
0.337447068854419	0.232315216341109\\
0.341679662777858	0.229679758995455\\
0.346064167222547	0.226260708384115\\
0.350006838735574	0.222591292917925\\
};
\addlegendentry{\trbdf}

\end{axis}
\end{tikzpicture}%

%% file: figures/high_order_traj_errors.tex
%
%
\begin{tikzpicture}

\begin{axis}[%
width=0.951\fwidth,
height=0.75\fwidth,
at={(0\fwidth,0\fwidth)},
scale only axis,
xmode=log,
xmin=0.007,
xmax=0.35,
xminorticks=true,
xlabel style={font=\color{white!15!black}},
xlabel={$h$},
ymode=log,
ymin=1e-06,
ymax=1000,
yminorticks=true,
axis background/.style={fill=white},
legend style={at={(0.03,0.5)}, anchor=west, legend cell align=left, align=left, draw=white!15!black}
]
\addplot [color=black, dotted]
  table[row sep=crcr]{%
0.0146861153387898	9.98169120109067e-07\\
0.0147590773878068	1.01783444363151e-06\\
0.0148903377550953	1.05587929018764e-06\\
0.0151962558832523	1.14636749805209e-06\\
0.0154290024367662	1.21967540331719e-06\\
0.0156324192756907	1.28618705761889e-06\\
0.0259889920659287	9.92448622749851e-06\\
0.0272323015390445	1.19399320399182e-05\\
0.0336902671714106	2.74262369146713e-05\\
0.0422595900184096	6.48735177359869e-05\\
0.0449671349844138	8.16700659328757e-05\\
0.0498759737559277	0.000119143209377742\\
0.055787894198809	0.000177410985990241\\
0.0611427180389079	0.000243606498197632\\
0.066614879114677	0.000325239472687064\\
0.0728987967807981	0.000437259994679122\\
0.0796680767446883	0.000580997790008523\\
0.101132853130827	0.00124016009182943\\
0.105621608281067	0.00144000628374307\\
0.109590381542605	0.00164806238133825\\
0.113237368472422	0.00187359073516452\\
0.116780701263661	0.00213352304151077\\
0.120285511043437	0.0024410635871552\\
0.123807693258545	0.00281374436090535\\
0.127446516841129	0.00328228937116418\\
0.131312753696466	0.00389544394120394\\
0.135504814423899	0.004726515437349\\
0.140229920571382	0.00592097787020218\\
0.145739743797227	0.00774651244142087\\
0.1527468939937	0.0109291784549204\\
0.164513948788841	0.0191966818058578\\
0.176676402932599	0.0327930301446338\\
0.183840299568759	0.0434689344763315\\
0.189547383427918	0.0531711883474632\\
0.194305264836055	0.0617204428909569\\
0.198316893013375	0.0688870057379727\\
0.201712634123019	0.0745932348502566\\
0.204586359050506	0.0788948477271774\\
0.206998810729389	0.0819300854195111\\
0.209012943137933	0.0839190010250836\\
0.210673424872781	0.0850917342281516\\
0.212034735660912	0.0856795608096768\\
0.21317008607798	0.0858755414097258\\
0.214190315531816	0.0857970945800076\\
0.215208963691094	0.085454241118439\\
0.216302446492949	0.0847634914889521\\
0.217495425640203	0.0835900912409192\\
0.218791245610056	0.0817645600158382\\
0.220173486107329	0.0791083201525955\\
0.221611708714808	0.0754635593171404\\
0.223075239744296	0.070694641166066\\
0.2245250790776	0.0647391786560849\\
0.225924580674961	0.0576120763782432\\
0.2272529923652	0.0493642701444875\\
0.228545742697736	0.0399023082139436\\
0.22992571797117	0.0316082470206646\\
0.23006396547091	0.0318753241816433\\
0.230273215070444	0.0339870019506082\\
0.230436781947714	0.0429614349988221\\
0.230198284358458	0.0618325897433514\\
0.229494919561536	0.0850762327512309\\
0.228390819058748	0.111848230148909\\
0.226896507188847	0.143125854353619\\
0.225003467315991	0.180312957138918\\
0.222704855563246	0.225011391117756\\
0.220007436929506	0.278857436763107\\
0.216943766092614	0.34316371575311\\
0.21356725239899	0.418739016531478\\
0.209937905656421	0.505924024740961\\
0.206099266360072	0.605067336263641\\
0.202075013643649	0.716811472655664\\
0.1978616565542	0.84265725420625\\
0.193422508601404	0.985713772387487\\
0.188676493791088	1.15201996125156\\
0.183424436283393	1.35548766049917\\
0.17675599808541	1.65433576228722\\
0.170259238994702	2.02795024859795\\
0.168373632136643	2.18147780274398\\
0.167534017732912	2.2861972551848\\
0.167335458274512	2.35904863603238\\
0.167567024492689	2.40831505965177\\
0.168081384419372	2.43975518544716\\
0.168770911748595	2.458072244062\\
0.169577845760298	2.46704226623075\\
0.170517855903801	2.46897923171004\\
0.171666277787647	2.4640753725344\\
0.173118647914027	2.45070004220377\\
0.174966062714026	2.42622055772293\\
0.177282266282462	2.38782445197302\\
0.180110983802209	2.33332464569407\\
0.183460515425668	2.26176096508058\\
0.187318105372201	2.17338663815237\\
0.191679748575013	2.06906913397403\\
0.196623026401101	1.94863697658479\\
0.202494968297064	1.80712319228915\\
0.221651513152799	1.42634374863039\\
0.225273073439767	1.37990730777114\\
0.228352675919736	1.34887114459474\\
0.231076812835629	1.32808245493168\\
0.233549819188479	1.31459583725423\\
0.235849983707906	1.3065097409568\\
0.238054050628329	1.30258888915802\\
0.24024947941216	1.30215408041779\\
0.242537357760999	1.30506743534246\\
0.245027513444332	1.31171436330087\\
0.247856149206302	1.32303582629691\\
0.251215780092984	1.34070568724538\\
0.255487161662681	1.36802465746361\\
0.261922295442961	1.4151046515786\\
0.276558587037306	1.52483930556988\\
0.283751301075902	1.57267515651817\\
0.290520996699776	1.6122184486008\\
0.297213121535052	1.64590377105663\\
0.30399957828854	1.67479207125047\\
0.31099357940635	1.69944866386308\\
0.318309158307992	1.72027132929031\\
0.326048445880089	1.73747125647766\\
0.334337159745095	1.75118266265927\\
0.343300273292168	1.76141014764448\\
0.350004043554996	1.76642924922691\\
};
\addlegendentry{\RadThree}

\addplot [color=blue, dashed]
  table[row sep=crcr]{%
0.0640933443727451	9.9967444570066e-07\\
0.0644506452018139	1.03658770897748e-06\\
0.0770038800126909	3.32359547164196e-06\\
0.0779584534115233	3.60641316252869e-06\\
0.079211948508372	4.01139829979237e-06\\
0.0883648882408492	8.42379234004009e-06\\
0.0898410253505953	9.45119702167366e-06\\
0.1030396216476	2.55550870887265e-05\\
0.116498451576278	6.81065121488735e-05\\
0.122304078794465	0.000103563482346348\\
0.128081071851009	0.000157021899950343\\
0.134243304610761	0.000244634465069355\\
0.140319439478475	0.000378533327738623\\
0.147866941791779	0.000650160941321467\\
0.155427681677305	0.00111505730717893\\
0.163695554360488	0.0020023004479746\\
0.173304086943932	0.00391523913350262\\
0.187977087447007	0.0104916101755132\\
0.198590422625616	0.0201867038987862\\
0.204126596239437	0.0273147393485723\\
0.207968106829272	0.0327702532862508\\
0.210738101440531	0.0365428914265687\\
0.212720907321368	0.0388108899595624\\
0.214098046295463	0.0399322393070335\\
0.215009914184303	0.0403312702966716\\
0.215583843280001	0.0403780407771211\\
0.21608800440191	0.0402442425273098\\
0.216648615410883	0.0398418153313166\\
0.217272327762679	0.0389549998869063\\
0.217911696736191	0.0372943349277571\\
0.218499674312196	0.0345191893274423\\
0.218964290994893	0.0302398332000461\\
0.219229116979941	0.0239957768455411\\
0.219175722526802	0.0154397312718457\\
0.21892265749136	0.012414252335927\\
0.218838039895649	0.0126424069966494\\
0.218541169276771	0.0158166552498054\\
0.217438542779835	0.0346170665168663\\
0.216320778743459	0.0539941103593464\\
0.214954337915589	0.0772661846012488\\
0.213322747843125	0.105005368834483\\
0.211424355209467	0.137711367757853\\
0.209264183053431	0.175847548795925\\
0.206854136199154	0.219773985058308\\
0.204210487418417	0.269740966787575\\
0.201352078412489	0.325887074245564\\
0.198297521575543	0.388278086834026\\
0.195064107616022	0.456934137799522\\
0.191665736939604	0.531893066542057\\
0.188113497140336	0.613227242957243\\
0.184416900620789	0.701043239369513\\
0.18058197128863	0.795553865330289\\
0.176613824522577	0.897043120510949\\
0.172515485791935	1.00591842161245\\
0.16828979881353	1.1226838188182\\
0.163938187242746	1.2479952631145\\
0.159460446605866	1.38269629481356\\
0.154858300468875	1.52774063237904\\
0.150132439150979	1.68430983210861\\
0.145283914412659	1.85381073014059\\
0.140315512552404	2.03786906759925\\
0.135229618341815	2.23845741462171\\
0.130030701564506	2.45786727355349\\
0.1247236143611	2.69885287517338\\
0.119316014643075	2.964618552855\\
0.113816759794556	3.25899754864367\\
0.108237164159871	3.58652538904129\\
0.102590071710804	3.95265697026115\\
0.0968911721281401	4.36389346291284\\
0.0911585320405425	4.82805440325744\\
0.0854113913492094	5.35470576031187\\
0.0796713687393415	5.95549935133158\\
0.0763639847711613	6.3423009452879\\
};
\addlegendentry{\RadFive}

\addplot [color=blue, dashed, forget plot]
  table[row sep=crcr]{%
0.0763639847711613	6.34230094528789\\
0.0686529272235586	7.38800709693616\\
0.0630676753156335	8.30467625415849\\
0.0575813587018133	9.37824672124948\\
0.0522197660989743	10.6457941550916\\
0.0470020773360151	12.1579655726529\\
0.0419497407010853	13.9824797657472\\
0.037083046600807	16.2129918679315\\
0.0324212580829643	18.9822464914961\\
0.0279814148190087	22.4851148036753\\
0.0237764280693236	27.0213910445586\\
0.019813041521664	33.0812856016247\\
0.0160865486586945	41.5428025922886\\
0.0125650789253304	54.2364014395402\\
0.0091533008078342	76.0784323126481\\
0.00699998655489947	101.137238604596\\
};
\addplot [color=red, dashdotted]
  table[row sep=crcr]{%
0.0390482872394351	9.99183931015124e-07\\
0.0440391273285853	1.74225094426779e-06\\
0.0452075867649297	1.96488952295226e-06\\
0.0518126368058897	3.684749721833e-06\\
0.0613141471109393	8.07566145299909e-06\\
0.0663144675989186	1.17376647829497e-05\\
0.0706421025822164	1.59826354572259e-05\\
0.0781475459768056	2.66768058042879e-05\\
0.0836429673343985	3.83457669031862e-05\\
0.0900269732866671	5.79826202033325e-05\\
0.0959003160807873	8.44029639534283e-05\\
0.100135519804641	0.000110402500659032\\
0.105135175566696	0.000151283676726638\\
0.110508044884289	0.000211767389089534\\
0.116686749398612	0.00031094535093036\\
0.122847411116763	0.000454889008417775\\
0.129219871409185	0.000672719297446328\\
0.135427812993071	0.000983339513840794\\
0.142051097240185	0.00147328783460958\\
0.14880112330732	0.00222547357720246\\
0.155487407001156	0.00335473047536491\\
0.162249843878463	0.00509793167856178\\
0.169012254384779	0.00778619104913264\\
0.175766099441235	0.0119666548810096\\
0.18250229647892	0.0185296391545666\\
0.189206464102957	0.028933777121031\\
0.19586755767605	0.0456155227238252\\
0.202442744046245	0.0725513951595585\\
0.208827630893471	0.115844563243053\\
0.214858187144873	0.184005228179811\\
0.22086793102485	0.299132114445625\\
0.225752197319977	0.436005833799733\\
0.227680319412269	0.486024336667804\\
0.229280306670824	0.517425244171493\\
0.230762697746376	0.537792013751325\\
0.232190767884442	0.550843912963553\\
0.233583147749606	0.558902350238713\\
0.234958092393274	0.563588413677817\\
0.236347561062975	0.565979258287788\\
0.237845016986577	0.566722198930213\\
0.239613670630319	0.565971761616669\\
0.241925949553972	0.563385562263037\\
0.245311529160396	0.557963974849982\\
0.250539783579318	0.548117371880432\\
0.255694957836302	0.537234197027266\\
0.259886627890878	0.526990791676201\\
0.263644069324129	0.516216057429618\\
0.267213444306821	0.504212747651011\\
0.270711535597933	0.490505706173695\\
0.274192533392705	0.474749569220304\\
0.277679498985329	0.456687632710759\\
0.281176633863715	0.43614627448972\\
0.284669099961875	0.413081225930296\\
0.288135002741861	0.387539298043637\\
0.291536124119172	0.359752324666994\\
0.294839411389552	0.32999262563257\\
0.29800784645082	0.298642430762869\\
0.301011447450586	0.266091416674044\\
0.303817013444164	0.232842373347529\\
0.306405282120301	0.199321063698583\\
0.308754180266199	0.166080644352483\\
0.310852681447979	0.133649134538936\\
0.312698211112288	0.102586728356802\\
0.314301757576088	0.0735012690839923\\
0.315728076712244	0.046812658200753\\
0.317279942562327	0.0277669695576628\\
0.317440027978705	0.0280940602070613\\
0.317809105145656	0.0312916863208155\\
0.319210062289483	0.0648936772266042\\
0.321318000232048	0.221754790924058\\
0.320819754476921	0.27668588173924\\
0.31977207684725	0.324510698111495\\
0.318265482806955	0.367432805431638\\
0.316338255286808	0.406872845885887\\
0.313974693195427	0.444520323896847\\
0.31105005803845	0.483256879834308\\
0.307120539163407	0.52997168975144\\
0.298899495057787	0.631716087899723\\
0.280911234078615	0.935289218536154\\
0.275264878020552	1.04410315588027\\
0.270063174951897	1.14385428030878\\
0.265010276100011	1.23846904160871\\
0.259946980922194	1.33024927072247\\
0.254697164278856	1.42233614950633\\
0.248934190715685	1.5211902106789\\
0.241554849463687	1.64922180727416\\
0.226058706734289	1.96700783455991\\
0.217953794813263	2.18547987049821\\
0.205100631625899	2.62642409776896\\
0.193724951503974	3.11127664002914\\
0.185543486297071	3.51162264891548\\
0.17809320669373	3.91331411955969\\
0.170950960499423	4.33278297768611\\
0.163948432510794	4.77896620609503\\
0.157000980779716	5.2590292458715\\
0.150063987307894	5.77977054806191\\
0.143114421264921	6.34834245074822\\
0.136140731881734	6.97285434060808\\
0.12914020506896	7.66264725353944\\
0.122117868920736	8.42854183564045\\
0.115082922568506	9.28348702909596\\
0.108048694958243	10.2430386384098\\
0.101030992774218	11.3262325654403\\
0.0940479451117479	12.5565666946775\\
0.0871193245205924	13.9634438674124\\
0.0802661402489503	15.584152417173\\
0.0735094949301689	17.4669237782513\\
0.0668708124960534	19.675073722103\\
0.0603709059275553	22.2936451322412\\
0.0540305594155692	25.4391193205635\\
0.0478708566779397	29.2748600239423\\
0.0419134161106278	34.0365033381632\\
0.0361815866551564	40.0742911352015\\
0.0307015820131145	47.927807152858\\
0.0274329235720058	54.0913381090066\\
};
\addlegendentry{\GauLegFour}

\addplot [color=red, dashdotted, forget plot]
  table[row sep=crcr]{%
0.0274329235720058	54.0913381090066\\
0.0223053147918916	67.3713223536492\\
0.0176618852846609	86.0112596394549\\
0.0134205335834974	114.264590811361\\
0.00964363686006197	160.292377340676\\
0.00699996172358029	222.02133269043\\
};
\addplot [color=black]
  table[row sep=crcr]{%
0.0988711765686779	9.99029348661864e-07\\
0.0992021519174039	1.02908815729325e-06\\
0.0999498064769086	1.09854484179144e-06\\
0.10021283997699	1.12438553326441e-06\\
0.100759415777247	1.18068712567072e-06\\
0.102549664794773	1.38212577388615e-06\\
0.102819123471979	1.41637965600569e-06\\
0.108888832901816	2.43164882363248e-06\\
0.109260817087764	2.51526171958866e-06\\
0.111848868639977	3.1726724597396e-06\\
0.112284214703106	3.30002308097314e-06\\
0.116831317504546	4.97339801523549e-06\\
0.117398130682313	5.23573813379475e-06\\
0.121614977906959	7.65905389547471e-06\\
0.126172516025004	1.15350765285689e-05\\
0.132134213815795	1.96157165800366e-05\\
0.137795103469989	3.22663636072184e-05\\
0.143633849186272	5.35305285816755e-05\\
0.150912283038233	9.96487380589507e-05\\
0.159756286747197	0.000209399530775531\\
0.167381043542226	0.000393714998348297\\
0.17660762896539	0.000837888135102362\\
0.186656214551195	0.00188658862624711\\
0.200635108629195	0.00564377577998572\\
0.211841851515665	0.0127072547018112\\
0.217883736472268	0.0187245735481028\\
0.222490578879377	0.0242645306349352\\
0.22623271659975	0.0290586140056654\\
0.229342129174534	0.0329122798026468\\
0.231911824949099	0.0357368992908433\\
0.234021993115918	0.0376206445557905\\
0.235722410782514	0.0387325856534288\\
0.237053110479167	0.0392834824290291\\
0.238063207706323	0.039481756912866\\
0.238858913272564	0.0394882069852801\\
0.239646532821512	0.0393527269501768\\
0.240540355058476	0.039013906603389\\
0.241577479382188	0.0383532108258617\\
0.242756235077447	0.0372213091985308\\
0.244051315362452	0.0354652157688967\\
0.245423966891182	0.0329605418981403\\
0.246830535669343	0.0296411536459897\\
0.248219517720256	0.0255485557930766\\
0.249539236336395	0.0208523917478964\\
0.250746183747751	0.0158336984729795\\
0.251806062094724	0.0108707381779206\\
0.252716197186835	0.00636086243158139\\
0.253631676699463	0.00328059177849441\\
0.253732453506313	0.00334071169895042\\
0.254019497968326	0.00410460457130295\\
0.255589783377958	0.0138839458818452\\
0.256998445019072	0.0251106731565251\\
0.258788843639771	0.0418388491644859\\
0.261062150885907	0.067525429658869\\
0.264128094743025	0.112667133488559\\
0.267574243770681	0.199997033220308\\
0.267784613192475	0.229177162288667\\
0.267347934302263	0.247296723479396\\
0.266543118290949	0.258338997587132\\
0.265545662779932	0.264640329491286\\
0.264489899094796	0.267829523900039\\
0.263470660574928	0.269098435393029\\
0.262473256436159	0.269251849937278\\
0.261336312061333	0.268518257980015\\
0.259892833799774	0.266589860339413\\
0.258015556910634	0.262883987136398\\
0.25565142354447	0.25681036899563\\
0.252871998095328	0.248110097776327\\
0.249861588055173	0.237036936511301\\
0.246784365293621	0.224008776305832\\
0.243714445946277	0.209270129029453\\
0.240639586080657	0.192794096282401\\
0.237468772247447	0.174261526016119\\
0.233874684687798	0.152295440800403\\
0.229319922673657	0.128047320517494\\
0.228769955235345	0.127194651968581\\
0.228464147517774	0.127669378223316\\
0.228114386051333	0.129838486368087\\
0.227721335535805	0.136541431504971\\
0.227334725487545	0.155065889091932\\
0.227035495381318	0.207618631336451\\
0.226950814789328	0.409082339533884\\
0.226656161160505	0.656115219840087\\
0.226002199941082	0.849974913531511\\
0.225017075967815	1.0344678015116\\
0.223707142370331	1.21847434307926\\
0.22206773181268	1.40734681125757\\
0.220085531311596	1.60584230473497\\
0.217739884031515	1.81899189518973\\
0.215007399734925	2.05221185245412\\
0.211870651061172	2.3109425815503\\
0.208330640060466	2.59978426956316\\
0.204415861576056	2.92155454834168\\
0.200178940194688	3.27700811934928\\
0.195682609191768	3.66541981959221\\
0.190981287308184	4.08592088757556\\
0.186113805585364	4.53833446999612\\
0.181104015107091	5.02352895190496\\
0.17596336758038	5.54360860805924\\
0.170697278114922	6.10165300205223\\
0.16530782391325	6.70171132553999\\
0.159795519689492	7.34884659860457\\
0.154159527528832	8.04935380713847\\
0.148400883958627	8.81061707067723\\
0.142521538094479	9.64148125163149\\
0.136525615907942	10.5523934480196\\
0.130419717305592	11.555711043138\\
0.124211666669421	12.6663668994197\\
0.117913056218304	13.9020099178355\\
0.111538039565984	15.2838778336737\\
0.105103394035109	16.8376752612418\\
0.0986291633348082	18.5945342358298\\
0.0921375308331214	20.5927438050299\\
0.0856533758354034	22.8795451445134\\
0.0792029285979821	25.5142206477319\\
0.072813685587603	28.5719339154028\\
0.0665140689770115	32.149221675335\\
0.0603323911015295	36.3723835088594\\
0.0542968096007152	41.409317122532\\
0.0522687156249924	43.3596527021617\\
};
\addlegendentry{\GauLegSix}

\addplot [color=black, forget plot]
  table[row sep=crcr]{%
0.0522687156249924	43.3596527021617\\
0.0461364164548579	50.2882942476027\\
0.0405646197516227	58.3825324962707\\
0.0352317034163197	68.50919517591\\
0.0301656264176881	81.4254955466678\\
0.0253948239257765	98.2790674877623\\
0.020948322486992	120.875437985547\\
0.0168564583983849	152.177989410909\\
0.0131500841866604	197.319840089311\\
0.00985949844492893	265.812328153015\\
0.00701269604958531	376.898548633045\\
0.00699998996398456	377.596839425714\\
};
\end{axis}
\end{tikzpicture}%

%% file: 0section-4-discussion.tex
\section{Discussion}
\label{sec:discussion}

Bifurcations, as a function of the timestep, are to be expected in the
solutions of the update equations defined by implicit numerical
methods for nonlinear IVPs. As demonstrated, such bifurcations can
occur in purely mathematical examples, such as the backward Euler
method applied to the differential equation $\dot q=q^2$, as well as
in more physically interesting models, such as the double pendulum.

For simulations where the accuracy of single trajectories is important
and variable step sizes are used, standard step-size control may
reduce the timestep to help prevent convergence to inconsistent
solutions. Such convergence, however, may persist even at small
timesteps, depending on how the iterative solver for a given method is
initialized.  Accordingly, convergence of an implicit method solver
alone is unreliable for assessing whether or not a solution is
consistent or whether a simulation should be carried on. In fact, it may
be possible to identify a critical timestep size $h_c$, corresponding
to the first point at which a bifurcation of the principal solution
branch occurs, such that timesteps $h > h_c$ can be considered to be
invalid and solutions obtained are by definition inconsistent.

In other instances, such as the long-time integration of symplectic
systems, where timesteps are often constant or more generally where
consideration of backward error is a driving motivation, specific
trajectories may not be as important as the general behaviour of the
system.  Even though the $\mathcal{O}(1)$ relative errors we observe
in trajectories at fold points in the double pendulum example imply
that the solutions may be inconsistent for timesteps smaller than the
critical timestep, it is helpful to be aware of the presence of
fundamental method- and state-dependent limitations on the size of the
valid timesteps. Simulations with timestep sizes that are larger than
the critical timestep cannot be expected to generate reasonable and
robust results in practice.

We used pseudo-arclength continuation to follow the solution branches
of the numerical solution to the double pendulum problem using several
common implicit numerical methods. Pseudo-arclength continuation can
be used in this fashion as a post hoc solution validation
technique. That is, given a simulation, we can use the
pseudo-arclength continuation procedure described to verify that the
solutions at each timestep lie on the principal solution branch. In
principle, this can be done in parallel with the computations of the
time-advanced state. It would be interesting to develop such a tool to
monitor timesteps of implicit methods and to flag if an inconsistent
solution is generated.

